\documentclass{amsart}
\usepackage{graphicx}
\usepackage{amsmath}

\newtheorem{thm}{Theorem}[section]
\newtheorem{cor}[thm]{Corollary}
\newtheorem{lem}[thm]{Lemma}
\newtheorem{prop}[thm]{Proposition}
\theoremstyle{definition}
\newtheorem{defn}[thm]{Definition}
\theoremstyle{definition}
\newtheorem{rem}[thm]{Remark}
\theoremstyle{definition}
\newtheorem{eg}[thm]{Example}
\numberwithin{equation}{section}

\newcommand{\Real}{\mathbb R}
\newcommand{\eps}{\varepsilon}
\newcommand{\lll}{\lambda}

\newcommand{\oct}{\mathbb O}

\newcommand{\J}{\mathcal{J}}
\newcommand{\jac}{\J_{(a,b)}}

\begin{document}

\title{(Semi-)Riemannian geometry of  (para-)octonionic projective planes}%
\author{Rowena Held, Iva Stavrov, Brian VanKoten}

\thanks{Special thanks to the John S. Rogers Summer Research Program at Lewis \& Clark College.}

\address{Department of Mathematical Sciences, Lewis \& Clark College, Portland, OR 97219, USA. Phone: 1-503-768-7560,\ Fax: 1-503-768-7668.}
\email{rmheld@lclark.edu, istavrov@lclark.edu, vankoten@lclark.edu}

\begin{abstract}
We use reduced homogeneous coordinates to study Riemannian geometry of the octonionic (or Cayley) projective plane. Our method extends to the para-octonionic (or split octonionic) projective plane, the octonionic projective plane of indefinite signature, and the hyperbolic dual of the octonionic projective plane;  we discuss these manifolds in the later sections of the paper. 
\end{abstract}
\maketitle
\section{Introduction}

I. Porteous \cite{portus} and H. Aslaksen \cite{aslaksen} have coordinatized the octonionic (Cayley) projective plane, $\oct P^2$, using the {\it reduced homogeneous coordinates}. Points in this model of $\oct P^2$ take the form $[u,v,w]$, with at least one of the octonions $u,v,w$ equal to $1$. For the appropriate coordinatization of lines in $\oct P^2$ the reader is referred to \cite{aslaksen}, where it is also shown that this model of projective plane geometry is non-Desarguesian. 

D. Allcock \cite{identification} provided an identification between the model of $\oct P^2$ involving reduced homogeneous coordinates and H. Freudenthal's \cite{frojdental} classic approach via Jordan algebras. Classically, the octonionic projective plane can be seen as a $16$-dimensional quotient manifold $F_4/Spin(9)$. This manifold can be equipped with a Riemannian metric with respect to which $\oct P^2$ is a $2$-point homogenous space.

In this paper we use the reduced homogeneous coordinates of Porteous to study the Riemannian geometry of $\oct P^2$. We explicitly write down the metric in terms of the reduced homogeneous coordinates, verify the homogeneity, and compute its curvature. We also provide an elementary approach to the other manifolds 
tightly related to the $\oct P^2$: the idefinite octonionic projective plane $\oct P^{(1,1)}$, the octonionic hyperbolic plane $\oct H^2$ and the para-octonionic projective plane $\oct'P^2$. In the concluding section of our paper we show that our models of the octonionic plane geometries are isometric to the classical models involving exceptional Lie groups. We do not provide explicit isometries, but use curvature classification results of Garcia-Rio, Vazques-Lorenzo and D. Kupeli \cite{spanishguys} regarding semi-Riemannian special Osserman manifolds. 

\section{Octonions and Para-octonions}

By the well-known result of Hurwitz \cite{hurevic} there are only 4 normed division algebras:
the real numbers $\mathbb{R}$, the complex numbers $\mathbb{C}$, the quaternions $\mathbb{H}$ and the octonions (or Cayley numbers) $\oct$. The octonions $\oct$ are an $8$-dimensional algebra which can be obtained from the quaternions using the Cayley-Dickson construction. That is, $\oct$ can be seen as the algebra $\mathbb{H}\oplus\mathbb{H}$ with multiplication given by 
\begin{equation}\label{oct}
(q_1,q_2)*(p_1,p_2)=(q_1p_1-\bar{p_2}q_2, p_2q_1+q_2\bar{p_1}).
\end{equation}
This algebra is not associative, as can be seen from the following:
\begin{alignat*}{5}
\big[(i,0)*(j,0)\big]*(0,1)&=(k,0)*(0,1)=(0,k),\\
(i,0)*\big[(j,0)*(0,1)\big]&=(i,0)*(0,j)=(0,-k).
\end{alignat*}
To measure the non-associativity of any three elements we can use the associator $[x,y,z]=x(yz)-(xy)z$.

The octonions have the property that any two elements generate an associative subalgebra \cite{Lie Algebras}. As a consequence, the associator is alternative, i.e. 
$$[a,b,c]=-[b,a,c]=-[a,c,b]=-[c,b,a].$$
Since $\bar{a}=2Re[a]-a$ and $[1,b,c]=0$ we have $[\bar{a},b,c]=-[a,b,c]$ and consequently 
$$\overline{[a,b,c]}=-[\bar{c},\bar{b},\bar{a}]=[c,b,a]=-[a,b,c].$$
In other words, the associator is always pure imaginary. As a consequence, the expression $Re[abc]$ well-defined, even though the expression $abc$ is not. It can also be shown that the expressions of the form $ab-ba$ are always pure imaginary.  Therefore, 
$$Re[abc]=Re[bca]=Re[cab].$$

The inner-product and the norm on $\oct$ are defined as
\begin{equation}\label{innerproduct}
<a,b>=\frac{a\bar{b}+b\bar{a}}{2}=Re[a\bar{b}]=Re[b\bar{a}],\ \ \ 
 <a,a>=|a|^2.
\end{equation}
This makes it clear that $|ab|^2=|a|^2|b|^2$. Moreover, we have
$$<ax,y>=<x,\bar{a}y> \text{\ and\ } <ax,ay>=|a|^2<x,y> \text{\ for all\ }a,x,y\in \oct.$$

We would like to be point out two identities we will use. The first one is 
 \begin{alignat}{5}
\label{re-stuff}
&Re[(ab)(cd)]+Re[(a\bar{c})(\bar{b}d)]=2Re[ad]Re[bc] \text{\ \ i.e.}\\
\label{<>stuff}
&<a\bar{b},c\bar{d}>+<a\bar{d},c\bar{b}>=2<a,c><b,d>.
\end{alignat}
This identity is basically a consequence of the alternativeness of the associator since
\begin{alignat*}{5}
&Re[(ab)(cd)]+Re[(a\bar{c})(\bar{b}d)]=Re[a\big(b(cd)\big)]+Re\big[a\big(\bar{c}(\bar{b}d)\big)\big]\\
=&Re[a(bc)d]+Re[a(\bar{c}\bar{b})d]+Re\big[a[b,c,d]\big]+
Re\big[a[\bar{c},\bar{b},d]\big]\\
=&Re[a(bc+\bar{c}\bar{b})d]=2Re[ad]Re[bc].
\end{alignat*}

The other identity we would like to point out is 
\begin{equation}\label{THE_ID}
(ab)(ca)=a(bc)a \text{\ \ for all\ } a,b,c\in \oct.
\end{equation}
The proof of this and many other identities involving the octonions can be found in \cite{frojdental}.

The para-octonions $\oct'$ can be constructed in a manner similar to (\ref{oct}). We define $\oct'$ to be the algebra  $\mathbb{H}\oplus\mathbb{H}$ with the  multiplication operation given by 
\begin{equation}\label{para-def}
(q_1,q_2)*(p_1,p_2)=(q_1p_1+\bar{p_2}q_2, p_2q_1+q_2\bar{p_1}). 
\end{equation}
This is a non-associative $8$-algebra whose  unit is $(1,0)$. Its standard basis vectors $(q,0)$ for $q\in\{i,j,k\}$ satisfy 
$(q,0)^2=1$, while $(0,w)$ for $w\in\{1,i,j,k\}$ satisfy $(0,w)^2=-1$. As with the octonions, this algebra has the property that any two elements generate an associative subalgebra.  This implies that the associator is alternative.

The inner-product on $\oct'$ can be defined using (\ref{innerproduct}). This inner-product is no longer positive definite but is of signature $(4,4)$. Indeed, the standard basis vectors of the type $(q,0)$ satisfy $|(q,0)|^2=1$ while the standard basis vectors of the type $(0,w)$ satisfy $|(0,w)|^2=-1$.  It is very important to notice that we still have
\begin{equation}\label{bigdeal}
|ab|^2=|a|^2|b|^2.
\end{equation}

Some identities in $\oct'$ which we will need are listed below. We do not know of a good reference in the literature for these, but they can all be easily derived from the definition of the multiplication (\ref{para-def}).

\begin{lem} \label{para-algebra}
If $a,b,c,x,y$ are arbitrary elements of $\oct '$ then 
\begin{enumerate}
\item $Re[ab]=Re[ba]$;
\item $Re\big{[}[a,b,c]\big{]}=0$ and 
$Re[abc]=Re[bca]=Re[cab]$ is well-defined;
\item $(ab)(ca)=a(bc)a$;
\item $<ax,y>=<x,\bar{a}y>$ and $<ax,ay>=|a|^2<x,y>$ with respect to the natural inner-product on $\oct'$.
\end{enumerate}
\end{lem}

\begin{proof}
We sketch the proof of property (3) in order to illustrate how the proofs of all of the other identities would go. We set $a=(a_1,a_2)$, $b=(b_1,b_2)$, $c=(c_1,c_2)$ and compute $(x',y'):=(ab)(ca)$ and $(x'',y''):=a(bc)a$ using the definition of the multiplication on $\oct'$. After some simplification we get
\begin{alignat*}{5}
&x'=a_1b_1c_1a_1+a_1\overline{c_2}b_2a_1+
2Re[b_1\overline{a_2}c_2]a_1+
2Re[\overline{b_2}a_2c_1]a_1+
|a_2|^2\overline{b_2}c_2+
|a_2|^2\overline{c_1}\overline{b_1},
\\ 
&x''=a_1b_1c_1a_1+a_1\overline{c_2}b_2a_1+2Re[\overline{b_1}\overline{c_2}a_2]a_1
+2Re[c_1\overline{b_2}a_2]a_1
+|a_2|^2\overline{c_1}\overline{b_1}
+|a_2|^2\overline{b_2}c_2.
\end{alignat*}
The corresponding terms are equal to one another due to symmetries of $Re$ such as  $Re[\alpha\beta\gamma]=Re[\gamma\alpha\beta]$ and $Re[\alpha\beta\gamma]=Re[\bar{\gamma}\bar{\beta}\bar{\alpha}]=Re[\bar{\beta}\bar{\alpha}\bar{\gamma}]=Re[\bar{\alpha}\bar{\gamma}\bar{\beta}]$. 
The equality between $y'$ and $y''$ can be shown in the same manner.
\end{proof}

\section{$\oct P^2$ as a Riemannian Manifold via Reduced Homogeneous Coordinates}\label{metricsection}

To describe the reduced homogeneous coordinates we consider a relation $\sim$ on $\mathbb{O}^3$; we say that $[a,b,c]\sim [d,e,f]$ if and only if there exists $\lll \in \mathbb{O}-\{0\}$ such that $a=d\lll, b=e\lll, c=f\lll$. This relation is symmetric and reflexive but due to non-associativity of octonions it is not necessarily transitive. To remedy this problem consider the following subsets of $\mathbb{O}^3$:
$$U_1=\{1\}\times\oct\times\oct,\ \ U_2=\oct\times\{1\}\times\oct,\ \ U_3=\oct\times\oct\times\{1\},$$
and their union $\mathcal{U}:=U_1\cup U_2\cup U_3$. 

\begin{lem}
The relation $\sim$ on $\mathcal{U}$ is an equivalence relation.
\end{lem}
\begin{proof} See \cite{aslaksen} or lemma \ref{equi}.\end{proof}
This lemma allows the following definition. 
\begin{defn}
The octonionic projective plane is the set of equivalence classes of $\mathcal{U}$ by the equivalence relation $\sim$.
$$\mathbb{O}\mathrm{P}^2=\mathcal{U}/_{\sim}$$
\end{defn}

\begin{thm}\label{mfld}
The octonionic projective plane $\mathbb{O} \mathrm{P}^2$ is a 16-dimensional simply connected manifold.
\end{thm}
\begin{proof}
We equip $\mathbb{O} \mathrm{P}^2$ with an atlas $(U_i/_{\sim}, \phi_i) (i=1,2,3)$, where the homeomorphisms $\phi_i$ are given by
\begin{alignat*}{5}
&\phi_1:U_1/_{\sim}\to\mathbb{R}^{16}\quad\phi_1([a,b,c])=(b,c);\\
&\phi_2:U_2/_{\sim}\to\mathbb{R}^{16}\quad\phi_2([a,b,c])=(a,c);\\
&\phi_3:U_3/_{\sim}\to\mathbb{R}^{16}\quad\phi_3([a,b,c])=(a,b).
\end{alignat*}
The transition functions $\phi_i\circ\phi_j^{-1}:\mathbb{R}^{16}\to\mathbb{R}^{16}$ 
are 
\begin{alignat*}{5}
&\phi_1\circ\phi_2^{-1}(a,b)=(a^{-1},ba^{-1})=\phi_2\circ\phi_1^{-1}(a,b);\\
&\phi_1\circ\phi_3^{-1}(a,b)=(ba^{-1},a^{-1})=\phi_3\circ\phi_1^{-1}(a,b);\\
&\phi_2\circ\phi_3^{-1}(a,b)=(b^{-1},ab^{-1})=\phi_3\circ\phi_2^{-1}(a,b),
\end{alignat*}
so $\mathbb{O} \mathrm{P}^2$ has a smooth 16-dimensional manifold structure. 

The open sets $U_1/_{\sim}, U_2/_{\sim}, U_3/_{\sim}$ are all simply connected since they are homeomorphic to $\Real^{16}$. The intersection $U_1/_{\sim}\cap U_2/_{\sim}=\big\{[1,u,v]\big| u\ne 0\big\}$ is simply connected as well since it is homeomorphic to $\big(\Real^8-\{0\}\big)\times \Real^8$, i.e. it is homotopy equivalent to $S^7\times \Real^8$. Thus, by Van Kampen theorem $U_1/_{\sim}\cup U_2/_{\sim}$ is simply connected. Furthermore, the intersection 
$\big(U_1/_{\sim}\cup U_2/_{\sim}\big)\cap U_3/_{\sim}=\big\{[u,v,1]\big| u,v\ne 0\big\}$ is simply connected because it is homotopy equivalent to $S^7\times S^7$. It follows from Van Kampen theorem that $\oct P^2=U_1/_{\sim}\cup U_2/_{\sim}\cup U_3/_{\sim}$ is also simply connected. 
\end{proof}

We now explain the Riemannian metric on $\mathbb{O}\mathrm{P}^2$. We first put a metric on each of the charts $U_1/_{\sim}$, $U_2/_{\sim}$, $U_3/_{\sim}$, and then check compatibility with respect to the transition maps.  

If $(u,v)$ are coordinate functions on these charts we set the metric as 
\begin{equation}\label{3.1}
ds^2=\frac{|du|^2(1+|v|^2)+|dv|^2(1+|u|^2)
-2Re[(u\bar{v})(dv d\bar{u})]}{(1+|u|^2+|v|^2)^2}\ .
\end{equation}
This choice of metric is motivated by the form of the Fubini-Study metric on complex and hyperbolic projective spaces; see \cite{Metric}. 


\begin{thm} \label{thm0}
The expression (\ref{3.1}) defines a Riemannian metric on $\mathbb{O}\mathrm{P}^2$.
\end{thm}

\begin{proof}
First observe that for all $(u,v)\in\oct^2$ the metric (\ref{3.1}) is positive definite. This is a consequence of the Cauchy-Schwartz inequality for the real inner product
$<x,y>=Re(x\bar{y})$ on $\oct$.  
Indeed, the inequality ensures that 
$$Re[(u\bar{v})(dv d\bar{u})]\le |u\bar{v}||dud\bar{v}|=|u||v||du||dv|$$
and so 
\begin{alignat*}{5}
ds^2&\ge \frac{|du|^2(1+|v|^2)+|dv|^2(1+|u|^2)-2|u||v||du||dv|}{(1+|u|^2+|v|^2)^2}\\
&=\frac{|du|^2+|dv|^2+(|du||v|-|dv||u|)^2}{(1+|u|^2+|v|^2)^2}>0.
\end{alignat*}

We now check that changes of coordinates preserve $ds^2$. Due to the symmetry of our transition functions (see the proof of Theorem \ref{mfld}) we can perform the calculation for the transition function $(u,v)=(x^{-1},zx^{-1})$. 
\begin{alignat*}{5}
ds^2
=&\Big{[}~|dx^{-1 }|^2(1+ |zx^{-1}|^2)+|d(zx^{-1})|^2(1+|x^{-1}|^2)\\
&\ \ \ -2Re\big[\big(x^{-1}\overline{x^{-1}}\bar{z}\big)\big(d(zx^{-1})d\overline{x^{-1}}\big)\big]\ \Big{]}\\
&~/(1+|x^{-1}|^2+|zx^{-1}|^2)^2.
\end{alignat*}
It follows from $d(xx^{-1})=dxx^{-1}+xdx^{-1}=0$ that $dx^{-1}=-x^{-1}dxx^{-1}$. Thus
$|dx^{-1}|^2=|x^{-1}|^4|dx|^2$ and
\begin{alignat}{5}
\label{messinxz1}
ds^2=&\Big{[}\ |x^{-1}|^4|dx|^2(1+|z|^2|x^{-1}|^2)+|dzx^{-1}+zdx^{-1}|^2(1+|x^{-1}|^2) \\
\label{messinxz2}
&\ \ -2|x^{-1}|^2Re\big[\bar{z}\big((dzx^{-1}+zdx^{-1})d\overline{x^{-1}}\big)\big]\ \Big{]} \\
&~/(1+|x^{-1}|^2+|z|^2|x^{-1}|^2)^2.
\end{alignat}
We apply $|a+b|^2=|a|^2+|b|^2+2Re[a\bar{b}]$ to see that 
\begin{alignat}{5}
|dzx^{-1}+zdx^{-1}|^2
=&|dzx^{-1}|^2+|zdx^{-1}|^2+2Re[(dz x^{-1})(d\overline{x^{-1}}\bar{z})]\\
=&|x^{-1}|^2|dz|^2+|x^{-1}|^4|z|^2|dx|^2+2Re[(dz x^{-1})(d\overline{x^{-1}}\bar{z})].\label{mess8/3}
\end{alignat}
The real part in equation (\ref{mess8/3}) can be simplified further. 
\begin{lem}\label{tech}
$Re[(dz x^{-1})(d\overline{x^{-1}}\bar{z})]= -|x^{-1}|^4 Re[(x\bar{z})(dzd\bar{x})].$
\end{lem}
\begin{proof} Since $Re[(ab)c]=Re[a(bc)]$ for all $a,b,c\in\oct$ we have
$$Re[(dz \bar{x})(\bar{x}^{-1}\bar{z})]=Re[dz(\bar{x}\bar{x}^{-1}z)]=Re[dz \bar{z}].$$
Taking the differential with respect to $x$ now produces
$$Re[(dzd\bar{x})(\bar{x}^{-1}\bar{z})]+Re[(dz\bar{x})(d\overline{x^{-1}}\bar{z})]=0.$$
Therefore, 
\begin{alignat*}{5}
Re[(dz x^{-1})(d\overline{x^{-1}}\bar{z})]
&=|x^{-1}|^2 Re[(dz\bar{x})(d\overline{x^{-1}}\bar{z})]=-|x^{-1}|^2 Re[(dzd\bar{x})(\bar{x}^{-1}\bar{z})]\\
&=-|x^{-1}|^4 Re[(dzd\bar{x})(x\bar{z})]=-|x^{-1}|^4 Re[(x\bar{z})(dzd\bar{x})];
\end{alignat*}
the last equality follows from $Re[ab]=Re[ba]$ for all $a,b\in\oct$.
\end{proof}
\noindent Using lemma \ref{tech} and identity (\ref{mess8/3}) we get  
$$|dzx^{-1}+zdx^{-1}|^2=|x^{-1}|^2|dz|^2+|x^{-1}|^4|z|^2|dx|^2-2|x^{-1}|^4 Re[(x\bar{z})(dzd\bar{x})].$$
We can also apply lemma \ref{tech} to simplify the line (\ref{messinxz2}). 
\begin{alignat*}{5}
&Re\big[\bar{z}\big((dzx^{-1}+zdx^{-1})d\overline{x^{-1}}\big)\big]
=Re\big[\bar{z}\big((dzx^{-1})d\overline{x^{-1}}\big)\big]+Re[\bar{z}zdx^{-1}d\overline{x^{-1}}]\\
=&Re\big[\big((dzx^{-1})d\overline{x^{-1}}\big)\bar{z}\big]+|z|^2|dx^{-1}|^2
=Re\big[(dzx^{-1})(d\overline{x^{-1}}\bar{z})\big]+|x^{-1}|^4|z|^2|dx|^2\\
=&-|x^{-1}|^4Re[(x\bar{z})(dzd\bar{x})]+|x^{-1}|^4|z|^2|dx|^2.
\end{alignat*}
Combining the last two expressions with (\ref{messinxz1}), (\ref{messinxz2}) and  $|x^{-1}|^2=|x|^{-2}$ we get
\begin{alignat*}{5}
ds^2
=&\frac{|dx|^2(1+|z|^2)+|dz|^2(1+|x|^2)-2Re[(x\bar{z})(dzd\bar{x})]}{(1+|x|^2+|z|^2)^2},
\end{alignat*}
which completes our proof.
\end{proof}

We proceed by discussing the components of the metric tensor $g$. Consider a point $P$ with coordinates $(u,v)$ and the coordinate frame $\{e_1,...,e_8,f_1,...,f_8\}$, where
$$e_i:=\partial_i,\ f_i:=\partial_{i+8},\ 1\le i\le 8.$$
It is immediate that 
$$g(e_i,e_j)=\delta_{ij}\frac{1+|v|^2}{(1+|u|^2+|v|^2)^2} \text{\ \ and\ \ }  
g(f_i, f_j)=\delta_{ij}\frac{1+|u|^2}{(1+|u|^2+|v|^2)^2}.$$ 
To describe $g(e_i, f_j)$ we employ $\{x_1,...,x_8\}$, the standard orthonormal basis of $\oct$. Since 
$$Re[(u\bar{v})(dvd\bar{u})]=<u\bar{v}, dud\bar{v}>=<(u\bar{v})dv,du>$$
we see that 
$$g(e_i,f_j)=g(f_j,e_i)=-\frac{<(u\bar{v})x_j,x_i>}{(1+|u|^2+|v|^2)^2}.$$

\section{$\oct P^2$ is homogeneous}\label{homogeneity}

In the absence of a convenient Riemannian submersion with $\oct P^2$ as a base space, we are forced to prove homogeneity of $\oct P^2$ {\it directly}. To be precise, we find a collection of isometries which act transitively on $\oct P^2$. Our isometries will be made out of local isometries described in the following proposition.

\begin{prop} \label{isometry1}
Let $r\in\mathbb{R}$ and let $\lll\in \oct$ be such that 
$r^2+|\lll|^2=1$. The map 
$\mathcal{R}_{r,\lll}:U_1/_{\sim} \to U_1/_{\sim}$ defined by
\begin{equation}\label{defiso}
[1,u,v]\mapsto [1,u',v'], \text{\ \ where\ \ }u'= ru+\lll v \text{\ \ and\ \ } v'= \bar{\lll}u-rv,
\end{equation}
is an isometry. 
\end{prop}

\begin{rem} One could define maps $\mathcal{R}_{r,\lll}$ on $U_2/_{\sim}$ and  $U_3/_{\sim}$ analogously; these would also  be isometries. The form of the maps $\mathcal{R}_{r,\lll}$ is motivated by the form of reflections in $\mathbb{H}^2$. Assuming the quaternionic inner-product on $\mathbb{H}^2$ is conjugate linear in the first entry, the reflection with respect to 
$(a,b)^\perp$ (where $|(a,b)|=1$)  takes the form of $\mathcal{R}_{r,\lll}$ with $r=|b|^2-|a|^2$ and $\lll=-2(a\bar{b})$. Note also that $\mathcal{R}_{r,\lll}^2=Id$. 
\end{rem}

\begin{proof}
The proof consists of a lengthy computation similar to the one performed in the proof of Theorem \ref{thm0}; we only point out those aspects of the computation which require non-trivial identities in $\oct$.  

Using $Re[u\bar{v}\bar{\lll}]=Re[u(\bar{v}\bar{\lll})]=
Re[(u\bar{v})\bar{\lll}]=
Re[\bar{\lll}(u\bar{v})]=Re[(\bar{\lll}u)\bar{v}]$ 
we easily obtain 
$$|u'|^2+|v'|^2=|u|^2+|v|^2 \text{\ \ and\ \ } |du'|^2+|dv'|^2=|du|^2+|dv|^2.$$
It now follows from the form of our metric (see (\ref{3.1})) that we only need to work with 
\begin{equation}\label{mess1001}
|du'|^2|v'|^2+|dv'|^2|u'|^2-2Re[(u'\overline{v'})(dv' d\overline{u'})].
\end{equation}
Direct substitution of (\ref{defiso}) and straightforward algebraic manipulation convert expression (\ref{mess1001}) to 
\begin{alignat*}{5}
&|du|^2|v|^2+|dv|^2|u|^2-2r^4Re[(u\bar{v})(dvd\bar{u})]\\
&+2r^2\Big\{Re\big[(u\bar{v})\big((\bar{\lll}du)(d\bar{v}\bar{\lll})\big)\big]+Re\big[(dvd\bar{u})\big((\lll v)(\bar{u}\lll)\big)\big]
-4Re[u\bar{v}\bar{\lll}]Re[dud\bar{v}\bar{\lll}]\Big\}\\
&-2Re\big[\big((\lll v)(\bar{u}\lll)\big)\big((\bar{\lll}du)(d\bar{v}\bar{\lll}\big)\big].
\end{alignat*}

Recall that $(ab)(ca)=a(bc)a$ for all $a,b,c\in\oct$ (see (\ref{THE_ID})). This means that 
$$Re\big[(u\bar{v})\big((\bar{\lll}du)(d\bar{v}\bar{\lll})\big)\big]=
Re\big[(u\bar{v})\big(\bar{\lll}(dud\bar{v})\bar{\lll})\big)\big]=
Re\big[\big((u\bar{v})\bar{\lll}\big)\big((dud\bar{v})\bar{\lll}\big)\big].$$
and similarly 
\begin{alignat*}{5}
Re\big[(dvd\bar{u})\big((\lll v)(\bar{u}\lll)\big)\big]=&
Re\big[(dvd\bar{u})\big(\lll (v\bar{u})\lll)\big)\big]\\
=&
Re\big[\big(\lll(dvd\bar{u})\big)\big(\lll (v\bar{u})\big)\big]=
Re\big[\big((u\bar{v})\bar{\lll}\big)\big((dud\bar{v})\bar{\lll}\big)\big].
\end{alignat*}
On the other hand, identity (\ref{re-stuff}) implies that 
$$Re\big[\big((u\bar{v})\bar{\lll}\big)\big((dud\bar{v})\bar{\lll}\big)\big]
+Re\big[\big((u\bar{v})(dvd\bar{u})\big)(\lll\bar{\lll})\big]=
2Re[u\bar{v}\bar{\lll}]Re[\bar{\lll}dud\bar{v}]$$
and consequently  
\begin{alignat*}{5}
&2r^2\big\{Re\big[(u\bar{v})\big((\bar{\lll}du)(d\bar{v}\bar{\lll})\big)\big]+Re\big[(dvd\bar{u})\big((\lll v)(\bar{u}\lll)\big)\big]
-4Re[u\bar{v}\bar{\lll}]Re[dud\bar{v}\bar{\lll}]\big\}\\
=&2r^2\big\{2Re\big[\big((u\bar{v})\bar{\lll}\big)\big((dud\bar{v})\bar{\lll}\big)\big]-4Re[u\bar{v}\bar{\lll}]Re[dud\bar{v}\bar{\lll}]\big\}=
-4r^2|\lll|^2Re\big[(u\bar{v})(dvd\bar{u})\big]. 
\end{alignat*}

In a similar fashion 
\begin{alignat*}{5}
&Re\big[\big((\lll v)(\bar{u}\lll)\big)\big((\bar{\lll}du)(d\bar{v}\bar{\lll}\big)\big]=Re\big[\big(\lll(v\bar{u})\lll\big)\big(\bar{\lll}(dud\bar{v})\bar{\lll}\big)\big]\\
=&Re\big[\big(\lll(v\bar{u})\big)\big(|\lll|^2(dud\bar{v})\bar{\lll}\big)\big]=|\lll|^2Re[\big((dud\bar{v})\bar{\lll}\big)\big(\lll(v\bar{u})\big)]\\
=&|\lll|^4Re[(dud\bar{v})(v\bar{u})]=|\lll|^4Re[(u\bar{v})(dvd\bar{u})].
\end{alignat*}

By combining these identities we get 
\begin{alignat*}{5}
&|du'|^2|v'|^2+|dv'|^2|u'|^2-2Re[(u'\overline{v'})(dv' d\overline{u'})]\\
=&|du|^2|v|^2+|dv|^2|u|^2-2(r^4+2r^2|\lll|^2+|\lll|^4)Re[(u\bar{v})(dvd\bar{u})]\\
=&|du|^2|v|^2+|dv|^2|u|^2-2Re[(u\bar{v})(dvd\bar{u})],
\end{alignat*}
which completes our proof that $\mathcal{R}_{r,\lll}$ is a (local) isometry.
\end{proof}

It is not a priori clear that the `` local reflections" $\mathcal{R}_{r,\lll}$ extend to globally defined maps on $\oct P^2$. We now show that they do.

\begin{prop} \label{extensions}
The local isometries $\mathcal{R}_{r,\lll}$ extend to (unique analytic) involutive isometries of $\oct P^2$.
\end{prop}

\begin{proof}
We first study possible extensions of $\mathcal{R}_{r,\lll}$ to $U_2/_{\sim}$. 
{\it Formally}, the ``reflection" $\mathcal{R}_{r,\lll}$ on $U_1/_{\sim}\cap U_2/_{\sim}$ looks like \begin{alignat*}{5}
[x,1,z]=&[1,x^{-1},zx^{-1}]\mapsto[1,rx^{-1}+\lll(zx^{-1}),\bar{\lll}x^{-1}-rzx^{-1}] \text{\ \ i.e.\ \ }\\
[x,1,z]\mapsto&\big[|x|^2\big(r\bar{x}+\lll(z\bar{x})\big)^{-1},1,\big(\bar{\lll}\bar{x}-rz\bar{x}\big)\big(r\bar{x}+\lll(z\bar{x})\big)^{-1}\big]\\
=&\big[|x|^2\big(\bar{\lll}\bar{x}-rz\bar{x}\big)^{-1}, \big(r\bar{x}+\lll(z\bar{x})\big)\big(\bar{\lll}\bar{x}-rz\bar{x}\big)^{-1}, 1\big].
\end{alignat*}
We continue by re-writing these expressions. Since 
$$\big(\bar{\lll}\bar{x}-rz\bar{x}\big)^{-1}=|x|^{-2}|\bar{\lll}-rz|^{-2}\big(x\lll-rx\bar{z}\big)$$ we have
$$
\big(r\bar{x}+\lll(z\bar{x})\big)\big(\bar{\lll}\bar{x}-rz\bar{x}\big)^{-1}
=\frac{r|x|^2\lll +\big(\lll(z\bar{x})\big)(x\lll)-r^2|x|^2\bar{z}-r|x|^2|z|^2\lll}
{|x|^2|\bar{\lll}-rz|^2}.
$$
Using (\ref{THE_ID}) we see that $\big(\lll(z\bar{x})\big)(x\lll)=
\lll (z\bar{x}x)\lll=|x|^2\lll z\lll$ and  so
$$
\big(r\bar{x}+\lll(z\bar{x})\big)\big(\bar{\lll}\bar{x}-rz\bar{x}\big)^{-1}=
|\bar{\lll}-rz|^{-2}\big\{r\lll+\lll z\lll-r^2\bar{z}-r|z|^2\lll\big\}.
$$
Furthermore,  
\begin{alignat*}{5}
|r\bar{x}+\lll(z\bar{x})|^2=&r^2|x|^2+2rRe[\bar{x}(x\bar{z})\bar{\lll}]+
|\lll|^2|z|^2|x|^2\\
=&|x|^2\{r^2+2rRe[\bar{z}\bar{\lll}]+|\lll|^2|z|^2\}=
|x|^2|r+\lll z|^2 \text{\ \ and \ \ }\\
\big(r\bar{x}+\lll(z\bar{x})\big)^{-1}=&|x|^{-2}|r+\lll z|^{-2}\big(rx+(x\bar{z})\bar{\lll}\big).
\end{alignat*}
A short computation in which we use  $(\bar{\lll}\bar{x})\big((x\bar{z})\bar{\lll}\big)=\bar{\lll}(\bar{x}x\bar{z})\bar{\lll}=|x|^2\bar{\lll}\bar{z}\bar{\lll}$ now yields
$$
\big(\bar{\lll}\bar{x}-rz\bar{x}\big)\big(r\bar{x}+\lll(z\bar{x})\big)^{-1}= |r+\lll z|^{-2}\big\{r\bar{\lll}+\bar{\lll}\bar{z}\bar{\lll}-r^2z-r|z|^2\bar{\lll}\big\}.
$$

Define
$$U_2':=\{[x,1,z]\in\oct P^2\big| r+\lll z\ne 0\},\ \ \ 
U_2'':=\{[x,1,z]\in\oct P^2 \big| \bar{\lll}-rz\ne 0\};$$ these are open subsets of 
$\oct P^2$. Since $\mathcal{R}_{r,\lll}$ is bijective we have 
$$r+\lll z\ne 0 \text{\ \ or\ \ } \bar{\lll} - rz\ne 0 \text{\ \ for all\  \ }z\in\oct.$$
Thus $U_2'\cup U_2''=U_2/_{\sim}$ and 
we may extend $\mathcal{R}_{r,\lll}$ to the whole of $U_2/_{\sim}$: 
\begin{equation}\label{ext1}
[x,1,z]\mapsto 
\begin{cases}
       \Big[\ \frac{rx+(x\bar{z})\bar{\lll}}{|r+\lll z|^2},\ 1,\ \frac{r\bar{\lll}+\bar{\lll}\bar{z}\bar{\lll}-r^2z-r|z|^2\bar{\lll}}{|r+\lll z|^2}\ \Big] 
       &  \text{\  for\  }[x,1,z]\in U_2', \cr
        &       \cr
         \Big[\ \frac{x\lll-rx\bar{z}}{|\bar{\lll}-rz|^2},\ \frac{r\lll+\lll z\lll-r^2\bar{z}-r|z|^2\lll}{|\bar{\lll}-rz|^2},\ 1\Big] & \text{\  for\  }[x,1,z]\in U_2'' .
 \end{cases}
\end{equation}
Similar extension process can be applied to $U_3/_{\sim}$. As the outcome we get
\begin{equation}\label{ext2}
[x,y,1]\mapsto 
\begin{cases}
       \Big[\ \frac{rx\bar{y}+x\bar{\lll}}{|ry+\lll|^2},\ 1,\ \frac{r|y|^2\bar{\lll}+\bar{\lll}y\bar{\lll}-r^2\bar{y}-r\bar{\lll}}{|ry+\lll|^2}\ \Big] 
       &  \text{\  for\  }[x,y,1]\in U_3', \cr
        &       \cr
         \Big[\ \frac{(x\bar{y})\lll-rx}{|\bar{\lll}y-r|^2},\ \frac{r|y|^2\lll+\lll \bar{y}\lll-r^2y-r\lll}{|\bar{\lll}y-r|^2},\ 1\Big] & \text{\  for\  }[x,y,1]\in U_3'' .
 \end{cases}
 \end{equation}
Here $U_3':=\{[x,y,1]\in\oct P^2\big| ry+\lll \ne 0\}$ and 
$U_3'':=\{[x,y,1]\in\oct P^2 \big| \bar{\lll}y-r\ne 0\}$; as above $U_3/_{\sim}=U_3'\cup U_3''$. 

The two extensions described above match on 
$U_1/_{\sim}\cap U_2/_{\sim}\cap U_3/_{\sim}$, 
which is open. Since all the coordinate expressions appearing in the extensions are rational, and since they match on the appropriate open subsets, the two extensions have to match on the entire $U_2/_{\sim}\cap U_3/_{\sim}$. Therefore, we have a well-defined, unique {\it analytic} extension of $\mathcal{R}_{r,\lll}$ to the entire $\oct P^2$. We shall use $\widetilde{\mathcal{R}_{r,\lll}}$ to denote this extension. 

We have $\mathcal{R}_{r,\lll}^2=Id$ on the open set $U_1/_{\sim}$. Since 
$\widetilde{\mathcal{R}_{r,\lll}}$ is the analytic extension of $\mathcal{R}_{r,\lll}$, we see that $\widetilde{\mathcal{R}_{r,\lll}}$ is an involution of $\oct P^2$. 

The components of our metric tensor are rational in the coordinates arising from the charts $U_i/_{\sim}$ (see (\ref{3.1})). The components of the pullback of our metric tensor using $\widetilde{\mathcal{R}_{r,\lll}}$ will also be rational due to the rational nature of $\widetilde{\mathcal{R}_{r,\lll}}$. Since $\mathcal{R}_{r,\lll}$ is an isometry of $U_1/_{\sim}$, the identity principle  shows that the extension of $\widetilde{\mathcal{R}_{r,\lll}}$ is a global isometry.
\end{proof}

\begin{eg}\label{rotation}
We illustrate the extensions on an example. Consider $\mathcal{R}_{1,0}$ and $\mathcal{R}_{\cos t,\sin t}$ on $U_3/_{\sim}$. Each of these ``reflections" extends to the whole of $\oct P^2$ and we get a global ``rotation"
$\mathcal{I}_t=\widetilde{\mathcal{R}_{1,0}}\circ\widetilde{\mathcal{R}_{\cos t,\sin t}}$. 
The ``rotation" $\mathcal{I}_t$ acts on $[x,y,1]$ as
 $$[x,y,1]\mapsto[\cos t\ x+\sin t\ y, -\sin t\ x+\cos t\ y, 1],$$ while it takes the point $[1,u,v]$ with $\cos t+\sin t\ u\ne 0$  to
 \begin{alignat*}{5}
      & [1, (-\sin t\ v^{-1} + \cos t\ uv^{-1})(\cos t\ v^{-1} +\sin t\ uv^{-1})^{-1},
       (\cos t\ v^{-1} +\sin t\ uv^{-1})^{-1}] \\
    =& [1, (-\sin t + \cos t\ u)(\cos t +\sin t\ u)^{-1},
       v(\cos t +\sin t\ u)^{-1}].
 \end{alignat*}
Thus, $\mathcal{I}_t([1,0,0])=[1,-\tan t,0]$ for all $t\ne \frac{\pi}{2}$. 
\end{eg}

\begin{thm}\label{homog1}
$\oct P^2$ is homogeneous.
\end{thm}

\begin{proof}
Due to symmetry it is enough to show that each point of $U_1/_{\sim}$ can be taken by an isometry to $[1,1,1]$ or, rather, that $[1,0,0]$ can be taken to any point of the form $[1,a,b]\ne [1,0,0]$. 

Let $R:=\sqrt{|a|^2+|b|^2}$ and distinguish two cases. 

{\it Case 1,} when $b=0$. In this case consider the (global) isometry arising from 
$\mathcal{R}_{0,1}\circ\mathcal{R}_{0,\frac{\bar{a}}{|a|}}$, where the ``reflections" involved are maps on $U_1/_{\sim}$. We have $$[1,R,0]=[1,|a|,0]\mapsto [1,0,a] \mapsto [1,a,0]=[1,a,b].$$

{\it Case 2,} when $b\ne 0$. In this case consider the (global) isometry arising from the composition 
$$\mathcal{R}_{\frac{|b|}{R},-\frac{a\bar{b}}{R|b|}}\circ\mathcal{R}_{0,-\frac{\bar{b}}{|b|}},$$
where the ``reflections" involved are maps on $U_1/_{\sim}$. 
This composition is taking 
$$[1,R,0]\mapsto \Big[1,0,-\frac{Rb}{|b|}\Big]
\mapsto[1,a,b].$$

The two cases above show that we can always find a global isometry taking $[1,R,0]$ to $[1,a,b]$. By precomposing such an isometry with  the ``rotation" $\mathcal{I}_t$ of  example \ref{rotation} (choose $t$ such that $\tan t=-R$) we see that $[1,0,0]$ can be taken to any point of the form $[1,a,b]$ via a global isometry.  
\end{proof}

\section{Curvature of $\oct P^2$}
\label{curvsect}

In this section we compute the Riemann curvature tensor of $\oct P^2$ with respect to the metric described in section \ref{metricsection}. We do so by understanding the curvature tensor at one particular point $P_0$ whose coordinates are $(0,0)$; this is sufficient as $\oct P^2$ is homogeneous. Note that at $P_0$ our metric is Euclidean, i.e. $g|_{P_0}=I$. In the next lemma we discuss the first and the second jets of $g$ at $P_0$.

\begin{lem} \label{jets}
The first jets of $g$ vanish at $P_0$. The only possibly non-vanishing second jets of $g$ at $P_0$ are listed below. 
\begin{enumerate}
\item $e_je_jg(e_i,e_i)=f_jf_jg(f_i,f_i)=-4$;
\item $f_jf_jg(e_i,e_i)=e_je_jg(f_i,f_i)=-2$;
\item $e_lf_kg(e_i,f_j)=-<x_l\overline{x_k}, x_i\overline{x_j}>$.
\end{enumerate}
\end{lem}

\begin{proof}
Using Maclaurin series we see that 
$$\frac{1}{(1+|u|^2+|v|^2)^2}=1-2(|u|^2+|v|^2)+O\Big(\big(|u|^2+|v|^2\big)^2\Big).$$
Hence the following quadratic approximations around $(u,v)=(0,0)$: 
\begin{alignat}{5} 
\label{approx1}
\frac{1+|u|^2}{(1+|u|^2+|v|^2)^2}&\approx1-2(|u|^2+|v|^2)+|u|^2=1-|u|^2-2|v|^2\\ 
\label{approx2}
\frac{1+|v|^2}{(1+|u|^2+|v|^2)^2}&\approx1-2(|u|^2+|v|^2)+|v|^2=1-2|u|^2-|v|^2\\
\label{approx3}
-\frac{<(u\bar{v})x_j,x_i>}{(1+|u|^2+|v|^2)^2}&\approx - <(u\bar{v})x_j,x_i>.
\end{alignat}
These approximations are accurate to fourth degree at $(u,v)=(0,0)$. Note that none of the approximations above has any linear terms. This implies that the first jets of the metric vanish at $P_0$. The jets listed under (1) and (2) are easily observed from the approximations (\ref{approx1}) and (\ref{approx2}). Finally, 
$$e_lf_k[<(u\bar{v})x_j,x_i>]=<(x_l\overline{x_k})x_j,x_i>
=<x_l\overline{x_k}, x_i\overline{x_j}>$$
yields the last equality.
\end{proof}

The Christoffel symbols are linear in the first jets of the metric. Therefore, the Christoffel symbols vanish at $P_0$. This fact, along with $g|_{P_0}=I$, considerably simplifies the expression for the components of the curvature tensor at $P_0$. Indeed, we have
\begin{equation}
\label{r's}
R_{\alpha\beta\gamma\delta}=R_{\alpha\beta\gamma}^{\ \ \ \ \delta}=
\Gamma_{\alpha\gamma;\beta}^{\delta}-\Gamma_{\beta\gamma;\alpha}^{\delta}
=\frac{1}{2}\big[g_{\beta\gamma;\alpha\delta}+g_{\alpha\delta;\beta\gamma}
-g_{\alpha\gamma;\beta\delta}-g_{\beta\delta;\alpha\gamma}\big].
\end{equation}
Our curvature computations will be based upon this formula. 

\begin{thm} \label{BIGTHM} The only possibly non-vanishing components of the curvature tensor are listed below. 
\begin{enumerate}
\item $R(e_i,e_j,e_i,e_j)=-R(e_i,e_j,e_j,e_i)=4$;
\item $R(f_i,f_j,f_i,f_j)=-R(f_i,f_j,f_j,f_i)=4$;
\item $R(e_i,e_j,f_k,f_l)=R(f_k,f_l,e_i,e_j)=-<x_i\overline{x_l},x_j\overline{x_k}>+<x_j\overline{x_l},x_i\overline{x_k}>$;
\item $R(e_i,f_j,e_k,f_l)=R(f_i,e_j,f_k,e_l)=<x_i\overline{x_j},x_k\overline{x_l}>$ and $R(e_i,f_j,f_l,e_k)=R(f_i,e_j,e_l,f_k)=
-<x_i\overline{x_j},x_k\overline{x_l}>$.
\end{enumerate}
\end{thm}

\begin{rem}
Note that one can summarize (1) and (2) with
\begin{alignat*}{5}
&R(e_i,e_j,e_k,e_l)=R(f_i,f_j,f_k,f_l)\\
=&-4<x_i,x_l><x_j,x_k>+4<x_j,x_l><x_i,x_k>.
\end{alignat*}
\end{rem}

\begin{proof}
It follows from the formula (\ref{r's}) and the previous lemma that 
\begin{alignat*}{5}
R(e_i,e_j,e_k,e_l)=&\frac{1}{2}\big[e_ie_lg(e_j,e_k)+e_je_kg(e_i,e_l)-
-e_je_lg(e_i,e_k)-e_ie_kg(e_j,e_l)\big]\\
=&-4(\delta_{il}\delta_{jk}-\delta_{jl}\delta_{ik}).
\end{alignat*}
In order for this curvature component to be non-zero we either need $i=l, j=k$ or $i=k, j=l$. In these cases we get the result listed under (1). 

The computation of $R(f_i,f_j,f_k,f_l)$ is completely analogous and will be omitted. Note however that claim (2) follows from (1) since
\begin{equation}\label{iso1}
[1,u,v]\mapsto[1,v,u]
\end{equation}
is a local isometry of $\oct P^2$ whose differential at $[1,0,0]$ interchanges $e_i$ with $f_i$.  

It is immediate from lemma \ref{jets} that  
$$R(e_i,e_j,e_k,f_l)=\frac{1}{2}\big[f_le_ig(e_j,e_k)+e_je_kg(f_l,e_i)-
f_le_jg(e_i,e_k)-e_ie_kg(f_l,e_j)\big]=0.$$
Using curvature symmetries and local isometry (\ref{iso1}) we now see that there are no non-vanishing curvature components of types
$R(e_*,e_*,e_*,f_*),\ R(e_*,e_*,f_*,e_*),$ $R(e_*,f_*,e_*,e_*),\ R(f_*,e_*,e_*,e_*),\ R(f_*,f_*,f_*,e_*),\ R(f_*,f_*,e_*,f_*),\ R(f_*,e_*,f_*,f_*),$ $R(e_*,f_*,f_*,f_*).$

We use formula (\ref{r's}) to find the remaining curvature components: 
\begin{alignat*}{5}
R(e_i,e_j,f_k,f_l)=&\frac{1}{2}\big[e_if_lg(e_j,f_k)+e_jf_kg(e_i,f_l)
-e_jf_lg(e_i,f_k)-e_if_kg(e_j,f_l)\big]\\
=&-<x_i\overline{x_l},x_j\overline{x_k}>+<x_j\overline{x_l},x_i\overline{x_k}> 
\text{\ \ \ and\ }\\
R(e_i,f_j,e_k,f_l)=&\frac{1}{2}\big[e_if_lg(e_k,f_j)+e_kf_jg(e_i,f_l)
-f_jf_lg(e_i,e_k)-e_ie_kg(f_j,f_l)\big]\\
=&-<x_i\overline{x_l},x_k\overline{x_j}>+2\delta_{ik}\delta_{jl}
\\=&-<x_i\overline{x_l},x_k\overline{x_j}>+2<x_i,x_k><x_j,x_l>.
\end{alignat*}
We have $-<x_i\overline{x_l},x_k\overline{x_j}>+2<x_i,x_k><x_j,x_l>=<x_i\overline{x_j},x_k\overline{x_l}>$ by (\ref{<>stuff}). Hence
$$R(e_i,f_j,e_k,f_l)=<x_i\overline{x_j},x_k\overline{x_l}>,$$ as claimed in part (4). 
The proof of our theorem is now easily obtained  
using curvature symmetries and the local isometry (\ref{iso1}). 
\end{proof}

We are now able to give a component-free description of the Riemann curvature tensor of $\oct P^2$. To express our result most efficiently we will identify the tangent space at $P_0[1,0,0]$ with pairs of octonions $(a,b)$ as follows
$$(a,b)=(\textstyle{\Sigma} a_ix_i,\textstyle{\Sigma} b_ix_i)\leftrightarrow\textstyle{\Sigma}a_ie_i+\textstyle{\Sigma}b_if_i.$$

\begin{cor}\label{R}
The Riemann curvature tensor of $\oct P^2$ at the point $P_0[1,0,0]$ is given by
\begin{alignat*}{5}
R\big( (a,b),(c,d),(e,f),(g,h)\big)=&4<a,e><c,g>-4<c,e><a,g>\\+&4<b,f><d,h>-4<d,f><b,h>\\
-&<e\bar{d},g\bar{b}>+<e\bar{b},g\bar{d}>-<c\bar{f}, a\bar{h}>+<a\bar{f},c\bar{h}>\\
-&<a\bar{d}-c\bar{b},g\bar{f}-e\bar{h}>.
\end{alignat*}
\end{cor}

\begin{proof}
The proof follows from the previous theorem using $\mathbb{R}$-multilinearity of the curvature tensor, the inner-product, the multiplication and the conjugation of octonions. 
\end{proof}

\begin{rem}
R. Brown and A. Grey \cite{grej} computed the curvature tensor of $\oct P^2$ using invariants of $Spin(9)$. Their result matches with ours after we alter our identification of $T_{P_0}\oct P^2$ with $\oct^2$. More precisely, had we used 
$$(a,\bar{b})=(\textstyle{\Sigma} a_ix_i,\textstyle{\Sigma} b_i\overline{x_i})\leftrightarrow\textstyle{\Sigma}a_ie_i+\textstyle{\Sigma}b_if_i$$
our expression would completely match with the one in \cite{grej}.
\end{rem}

\begin{rem}
Given an algebra structure on $\mathbb{R}^n$ (e.g Cayley-Dickson algebras) the formula of corollary \ref{R} gives rise to an algebraic curvature tensor (see \cite{elbook}) on $\mathbb{R}^n$.
\end{rem}

\section{Para-octonionic projective plane $\oct'P^2$}\label{para}

As the first step of creating our restricted homogeneous coordinates we define a subset of $\oct'$
 $$\oct'_+:=\{x\in \oct' \big{|} |x|^2>0\}.$$
It follows from $|ab|^2=|a|^2|b|^2$ and 
$|\frac{\bar{a}}{|a|^2}|^2=\frac{\bar{a}}{|a|^2}\frac{a}{|a|^2}=\frac{1}{|a|^2}$ that $\oct'_+$ is closed under multiplication and inverses. From a topological point of view we have the homotopy equivalences 
\begin{equation}\label{homotopyeq} 
\oct'_+=\big\{(x_1,....,x_8)\in\Real^8\big| x_1^2+...+x_4^2-x_5^2-...-x_8^2>0\big\}\simeq \Real^4-\{0\} \simeq S^3.
\end{equation}
These equivalences can be justified using the linear homotopy $$H(x_1,...,x_8;t)=(x_1,...,x_4,tx_5,...,tx_8), t\in[0,1].$$
 
Consider the relation $\sim$ on $\oct'^3$ for which $[a,b,c]\sim [d,e,f]$ if and only if there exists $\lambda \in \oct'_+$ such that $a=d\lambda, b=e\lambda, c=f\lambda$. This relation is reflexive and symmetric but not necessarily transitive. Define 
\begin{alignat*}{5}
U_1:=&\{(1,y_1,z_1)\in\mathbb{O}'^3\big{|} 
1+|y_1|^2+|z_1|^2>0\},\\ 
U_2:=&\{(x_2,1,z_2)\in\mathbb{O}'^3\big{|} |x_2|^2+1+|z_2|^2>0\},\\  
U_3:=&\{(x_3,y_3,1)\in\mathbb{O}'^3\big{|} |x_3|^2+|y_3|^2+1>0\}
\end{alignat*}
and consider the restriction of $\sim$ to the union  
$\mathcal{U}:=U_1\cup U_2\cup U_3$. 

\begin{lem}\label{equi}
The relation $\sim$ on $\mathcal{U}$ is an equivalence relation.
\end{lem}
\begin{proof}
Suppose  $[1,y_1,z_1]\sim[x_2,1,z_2]$   and   $[x_2,1,z_2]\sim[x_3,y_3,1]$. This means that  for some $\lambda_1,\lambda_2 \in \oct'_+$ we have 
\begin{equation}\label{6eqns}
\begin{cases}
       1=x_2\lambda_1       &  \cr
       y_1=\lambda_1&       \cr
       z_1=z_2\lambda_1  &  
 \end{cases}
 \text{\ \ and\ \ }
 \begin{cases}
       x_2=x_3\lambda_2 &   \cr
       1=y_3\lambda_2 &       \cr
       z_2=\lambda_2 .
 \end{cases}
\end{equation}
Since $\lll_1\in\oct'_+$ we have 
$y_1=\lambda_1\in\mathbb{O}'_+$ and   $x_2=\lambda_1^{-1}\in\mathbb{O}'_+$.    Similarly, $y_3,z_2\in\mathbb{O}'_+$. We now see that $x_3=x_2\lambda_2^{-1}\in\mathbb{O}'_+$,  $z_1=z_2\lambda_1\in\mathbb{O}'_+$ and hence  all of the $x_i$, $y_i$, $z_i$ involved are invertible. By eliminating $\lll_i$ from  (\ref{6eqns}) we easily obtain 
\begin{align*}
&x_3=x_2z_2^{-1}=y_1^{-1}(z_1y_1^{-1})^{-1}=y_1^{-1}(y_1z_1^{-1})=z_1^{-1}\\
&y_3=z_2^{-1}=(z_1y_1^{-1})^{-1}=y_1z_1^{-1}.
\end{align*}
We now set $\lambda_3^:=x_3^{-1}$. Then $\lambda_3\in \mathbb{O}'_+$ and 
$$1=x_3\lambda_3,\quad y_1=y_3\lambda_3,\quad z_1=\lambda_3.$$
Therefore, $[1,y_1,z_1]\sim[x_3,y_3,1]$ and the relation $\sim$ is transitive on $\mathcal{U}$.
\end{proof}
We may now consider the equivalence classes:
\begin{defn}
The para-octonionic projective plane is the set of equivalence classes of $\mathcal{U}$ by the equivalence relation $\sim$.
$$\oct'P^2=\mathcal{U}/_{\sim}$$
\end{defn}
\begin{thm} The para-octonionic projective plane $\oct'P^2$ is a $16$-dimensional simply connected manifold.
\end{thm}
\begin{proof}
The charts for $\oct'P^2$ arise from the sets $U_i/_{\sim},\ i=1,2,3$; the situation is analogous to the case of $\oct P^2$ (see theorem \ref{mfld}) and the details will be omitted. We proceed by proving that $\oct'P^2$ is simply connected.

The linear homotopy $H\big((u,v),t\big)=(tu,tv),\ t\in[0,1]$ shows that the open sets $U_i/_{\sim}\approx\big\{|u|^2+|v|^2>-1\big\}\subset \Real^{16},\ i=1,2,3$ are contractible. Now consider 
$$U_1/_{\sim}\cap U_2/_{\sim}=\big\{[1,u,v]\big|1+|u|^2+|v|^2>0, |u|^2>0\big\}.$$
Using the linear homotopy $H\big((u,v),t\big)=(u,tv), t\in[0,1]$ we see that $U_1/_{\sim}\cap U_2/_{\sim}$ is homotopy equivalent to $\oct'_+$ and  $S^3$ (see (\ref{homotopyeq})). In other words, $U_1/_{\sim}\cap U_2/_{\sim}\simeq S^3$ and the set $U_1/_{\sim}\cap U_2/_{\sim}$ is simply connected. By Van Kampen theorem  $U_1/_{\sim}\cup U_2/_{\sim}$ is simply connected as well.
 
Next consider 
$$\big(U_1/_{\sim}\cup U_2/_{\sim}\big)\cap U_3/_{\sim}=\big\{[u,v,1]\big| |u|^2>0,\ |v|^2>0\big\}\approx\oct'_+\times\oct'_+.$$ 
Since $\oct'_+\times\oct'_+\simeq S^3\times S^3$, due to (\ref{homotopyeq}), the set  $\big(U_1/_{\sim}\cup U_2/_{\sim}\big)\cap U_3/_{\sim}$ is simply connected. The Van Kampen theorem now shows that 
$\oct'P^2=U_1/_{\sim}\cup U_2/_{\sim}\cup U_3/_{\sim}$ is also simply connected.
\end{proof}

The semi-Riemannian geometry of the para-octonionic projective plane $\oct'P^2$ is to a great extent analogous to the Riemannian geometry on  $\oct P^2$. We put a metric on $\oct'P^2$ by first defining it on each of the charts $U_1/_{\sim}$, $U_2/_{\sim}$, $U_3/_{\sim}$. To be precise, we set  
\begin{equation}\label{6.?}
ds^2=\frac{|du|^2(1+|v|^2)+|dv|^2(1+|u|^2)
-2Re[(u\bar{v})(dv d\bar{u})]}{(1+|u|^2+|v|^2)^2}\ ,
\end{equation}
where $(u,v)$ are coordinate functions on any of the three charts.
It is now necessary to verify that the metrics are non-degenerate and that on any chart intersections the two metrics are the same. 

Verifying the compatibility on the overlaps carries over from $\oct P^2$ without any modification. Therefore, we only need to study the issue of non-degeneracy. To do this most efficiently we consider the metric tensor components. We will use the standard orthonormal basis $\{x_1,...,x_8\}$ of $\oct'$ and the coordinate frame
$$\{e_1,...,e_8,f_1,...,f_8\} \text{\ \ where\ \ }
e_i:=\partial_i,\ f_i:=\partial_{i+8},\ 1\le i\le 8.$$ The basis  $\{x_1,...,x_8\}$ is ``orthonormal" in the sense that $<x_i,x_j>=\delta_{ij}\eps_i$, where
$$
\eps_i:=
\begin{cases}
       1 &  i\le 4 \cr
       -1 & i\ge 5.
\end{cases}
$$
Due to the indefiniteness of the inner product on $\oct'$ we have 
\begin{alignat*}{5}
g(e_i,e_j)=&\delta_{ij}\eps_i\frac{1+|v|^2}{(1+|u|^2+|v|^2)^2},\ \ \  
g(f_i,f_j)=\delta_{ij}\eps_i\frac{1+|u|^2}{(1+|u|^2+|v|^2)^2} \\
g(e_i,f_j)=&g(f_j,e_i)=-\frac{<(u\bar{v})x_j,x_i>}{(1+|u|^2+|v|^2)^2}.\ \ 
\end{alignat*}
We can now write the metric tensor $g$ in the matrix form 
$$\mathcal{M}=
\frac{1}{(1+|u|^2+|v|^2)^2}
\left[\begin{array}{cc}(1+|v|^2)G & A \\ A^T & (1+|u|^2)G\end{array}\right]$$
where $A_{ij}=-<(u\bar{v})x_j,x_i>$, and where $G$ is the $8\times 8$ matrix $\left(\begin{array}{cc}I & 0 \\0 & -I\end{array}\right)$. 
Note the following property of the matrix $A$. 
\begin{lem}\label{kinda_orth}
$AGA^T=|u|^2|v|^2G=A^TGA.$
\end{lem}
\begin{proof}
We know $A_{ij}=-<(u\bar{v})x_j,x_i>=-<(v\bar{u})x_i,x_j>$. Therefore, 
$$(AGA^T)_{ij}=\sum^8_{k=1}\eps_k<(v\bar{u})x_i,x_k><(v\bar{u})x_j,x_k>.$$
If we consider that in general $\sum_i\eps_i<w,e_i><v,e_i>=<w,v>$ we can see that 
$$(AGA^T)_{ij}=<(v\bar{u})x_i,(v\bar{u})x_j>.$$
By lemma \ref{para-algebra} we now have 
$(AGA^T)_{ij}=|u|^2|v|^2\eps_i\delta_{ij}$
and so $$AGA^T=|u|^2|v|^2G.$$ 
This identity implies that $A$ is invertible whenever $|u|^2|v|^2\neq0$ and that  
$$A^{-1}=\frac{1}{|u|^2|v|^2}GA^TG.$$ 
In particular, 
$\frac{1}{|u|^2|v|^2}GA^TGA=Id$
for $|u|^2|v|^2\neq 0$ and by continuity $$A^TGA=|u|^2|v|^2G$$ 
for all $u,v\in\oct'$.
\end{proof}
We are now ready to prove the non-degenracy and compute the signature of the metric given by (\ref{6.?}). 
\begin{prop}\label{signature}
The inner product defined by matrix $\mathcal{M}$ is of signature (8,8).
\end{prop}
\begin{proof}
We first show that $\mathcal{M}$ is non-degenerate. Suppose there exists 
$\vec{v}=(\vec{r},\vec{s})$ such that  $\mathcal{M}\vec{v}=0$ i.e.
$$\left[\begin{array}{cc}(1+|v|^2)G & A \\ A^T & (1+|u|^2)G\end{array}\right]\left[\begin{array}{cc}\vec{r}\\ \vec{s} \end{array}\right]=0.$$
Then we have the system of equations
$$\begin{cases}
(1+|v|^2)G\vec{r}+A\vec{s}=0 & \cr
A^T\vec{r}+(1+|u|^2)G\vec{s}=0 & 
\end{cases}
$$
which further implies 
\begin{equation*}
\begin{cases}
(1+|u|^2)(1+|v|^2)G\vec{r}-(AGA^T)\vec{r}=0 & \cr
(1+|u|^2)(1+|v|^2)G\vec{s}-(A^TGA)\vec{s}=0. &
\end{cases}
\end{equation*}
Lemma \ref{kinda_orth} converts these equations  into
$$(1+|u|^2+|v|^2)G\vec{r}=0 \text{\ \ and\ \ }
(1+|u|^2+|v|^2)G\vec{s}=0.$$
Thus $\vec{r}=\vec{s}=0$ i.e. $\vec{v}=0$ and our inner product is non-degenerate. 

The signature of a non-degenerate metric is the same at every point. Thus we only need to consider the signature at one point. For convenience let us consider the point with coordinates $u=v=0$; at this point $ds^2=|du|^2+|dv|^2$. Since the para-octonionic inner product used to find both $|du|^2$ and $|dv|^2$ is of signature (4,4), our metric must be of signature (8,8). 
\end{proof}

It can be proved that this metric makes $\oct'P^2$ a homogeneous manifold. The basic idea behind the proof is the same as for $\oct P^2$. One first studies the maps of the form (\ref{defiso}). The proof of proposition \ref{isometry1} carries over to $\oct'P^2$ as all the octonionic identities we use in the proof also hold in $\oct'$ (see lemma \ref{para-algebra}).
So, the maps (\ref{defiso}) are local isometries of $\oct'P^2$. The issue of extending these isometries to the entire $\oct'P^2$ has certain technical subtleties. For example, it is not obvious from (\ref{ext1})-(\ref{ext2}) that the indicated extensions are even well-defined. 

Note that for $[1,u',v']=\mathcal{R}_{r,\lll}[1,u,v]$ we have 
$1+|u'|^2+|v'|^2=|1|+|u|^2+|v|^2$. Consequently, if $[x',1,z']$ is the {\it formal} image of $[x,1,z]$ under (\ref{ext1}), then 
$$|x'|^2+1+|z'|^2=\frac{|x'|^2}{|x|^2}\big(|x|^2+1+|z|^2\big).$$
It follows from  $|x'|^2=|x|^2|r+\bar{z}\bar{\lll}|^2$
and the rational nature of the map (\ref{ext1}) that 
$$|x'|^2+1+|z'|^2=
\frac{1}{|r+\lll z|^2}\big(|x|^2+1+|z|^2\big) \text{\ \ for all\ \ } x,z\in\oct' \text{\ \ with\ \ } |r+\lll z|^2\ne0.$$
A similar argument shows that if $[x',y',1]$ is the formal image of $[x,1,z]$ under (\ref{ext1}) then 
$$|x'|^2+|y'|^2+1=\frac{1}{|\bar{\lll}-rz|^2}\big(|x|^2+1+|z|^2\big).$$
Therefore, to ensure well-definedness of (\ref{ext1}) on $\oct'P^2$ we need to restrict our attention to the sets 
$$U_2':=\{[x,1,z]\in \oct'P^2\big||r+\lll z|^2>0\} \text{\ \ and\ \ } U_2'':=\{[x,1,z]\in \oct'P^2\big||\bar{\lll}-rz|^2>0\}.$$ 
Of course, these sets 
may no longer cover $U_2/_{\sim}$. However, since $$|r+\lll z|^2+|\bar{\lll}-rz|^2=1+|z|^2 
\text{\ \ and\ \ }1+|x|^2+|z|^2>0$$ we must have 
$$[x,1,z]\not\in U_2'\cup U_2''\Longrightarrow [x,1,z]\in U_1/_{\sim} \text{\ \ i.e. \ \ } U_1/_{\sim}\cup U_2'\cup U_2''=U_1/_{\sim}\cup U_2/_{\sim}.$$
If we consider 
$$U_3':=\{[x,y,1]\in \oct'P^2\big||ry+\lll |^2>0\} \text{\ \ and\ \ } U_3'':=\{[x,y,1]\in \oct'P^2\big||\bar{\lll}y-r|^2>0\}$$ 
we have 
$$\oct'P^2=U_1/_{\sim}\cup U_2'\cup U_2''\cup U_3'\cup U_3''.$$
Using this set-up we can easily see that 
the proposition \ref{extensions} carries over to $\oct'P^2$. 
\begin{thm}
$\oct' P^2$ is homogeneous. 
\end{thm}
\begin{proof}
We follow the basic idea of the proof of theorem \ref{homog1}. For a point $[1,a,b]$ such that $|b|^2\ne 0$ and $|a|^2+|b|^2\ne 0$ consider
$$K=\begin{cases}
       \sqrt{|a|^2+|b|^2}
       &  \text{\  if\  }|a|^2+|b|^2>0, \cr
        &       \cr
         -\sqrt{-|a|^2-|b|^2}
       &  \text{\  if\  }|a|^2+|b|^2<0 
 \end{cases}
 \text{\ \ and\ \ }
 L=\begin{cases}
       \sqrt{|b|^2}
       &  \text{\  if\  }|b|^2>0, \cr
        &       \cr
         -\sqrt{-|b|^2}
       &  \text{\  if\  }|b|^2<0 ;
 \end{cases}$$
we have $L|L|=|b|^2$. The composition 
$$\mathcal{R}_{-\frac{L}{K},\frac{a\bar{b}}{|L|K}}\circ\mathcal{R}_{0,\frac{\bar{b}}{L}}$$
takes $[1,K,0]$ to $[1,a,b]$. Since $[1,0,0]$ can be taken to any of the points of the form $[1,K,0]$ by a ``rotation" (see example \ref{rotation}), we see that $[1,0,0]$ can be taken to any of the points $[1,a,b]$ with $|b|^2\ne 0$ and $|a|^2+|b|^2\ne 0$. 

On the other hand, any point $[1,a,b]$ with $|b|^2=0$ or $|a|^2+|b|^2=0$ can be taken to a point of the form $[1,x,y]$ with $|y|^2\ne 0$ and $|x|^2+|y|^2\ne 0$ using  ``rotations" of the example \ref{rotation}.
More precisely, let $t_1,t_2$ be real numbers such that the function 
$$t\mapsto \cot^2t+2Re[b]\cot t +|b|^2$$
is positive on the interval $[t_1,t_2]$. We have that  
$$|\cos t+\sin t\ b|^2=\cos^2t+2\sin t\cos tRe[b]+\sin^2t |b|^2>0 \text{\ \ for all\ \ } t\in[t_1,t_2].$$ Therefore, we may consider the ``rotation" which on the neighborhood of $[1,a,b]$ looks like
\begin{equation}\label{works?}
[1,u,v]\mapsto[1,u(\cos t+\sin t\ v)^{-1}, 
(-\sin t +\cos t\ v)(\cos t+\sin t\ v)^{-1}].
\end{equation}
In the case when $|b|^2\ne 1$ or $Re[b]\ne 0$ the expression
$$\cos^2t+2\sin t\cos tRe[b]+\sin^2t |b|^2=
\frac{1-|b|^2}{2}\cos (2t)+Re[b]\sin(2t)+\frac{1+|b|^2}{2}$$ is non-constant and there exists $t\in[t_1,t_2]$ for which  
$$|-\sin t +\cos t b|^2=1+|b|^2-\big(\cos^2t+2\sin t\cos tRe[b]+\sin^2t |b|^2
\big)\ne 0$$
and also 
\begin{alignat*}{5}
&|a|^2+|-\sin t +\cos t b|^2\\
=&1+|a|^2+|b|^2-\big(\cos^2t+2\sin t\cos tRe[b]+\sin^2t |b|^2
\big)\ne 0.
\end{alignat*}
In particular, if 
$|b|^2\ne 1$ or $Re[b]\ne 0$ there exists an ``angle " $t$ such that the global ``rotation" arising from (\ref{works?}) takes 
the point $[1,a,b]$ to a point $[1,x,y]$ for which $|y|^2\ne 0$ and $|x|^2+|y|^2\ne 0$. The case when $|b|^2=1$ (and consequently  $|a|^2=-1$) can be handled by precomposing the ``rotation" (\ref{works?}) with a ``reflection" $\mathcal{R}_{\cos s, \sin s}$ on $U_1/_{\sim}$ where 
$|\sin s\ a-\cos s\ b|^2=\cos (2s) -Re[a\bar{b}]\sin (2s)\ne 1$.
\end{proof}

The curvature of $\oct ' P^2$ can be studied using methods of section \ref{curvsect}. As $\oct' P^2$ is homogeneous we may restrict our attention to the point $P_0[1,0,0]$ whose coordinates are $(0,0)$. The curvature computation is fairly easy at this point: our metric at $P_0$ is pseudo-Euclidean (that is,  $g|_{P_0}=G$), the first jets of the metric vanish and the only possibly non-vanishing second jets are: 
\begin{itemize}
\item $e_je_jg(e_i,e_i)=f_jf_jg(f_i,f_i)=-4\eps_i\eps_j$;
\item $f_jf_jg(e_i,e_i)=e_je_jg(f_i,f_i)=-2\eps_i\eps_j$;
\item $e_lf_kg(e_i,f_j)=-<x_l\overline{x_k}, x_i\overline{x_j}>$.
\end{itemize}
These equalities can be verified in the same manner as in lemma \ref{jets}.

The Christoffel symbols vanish at $P_0$ because they are linear in the first jets of the metric. Hence the following expression for the components of the curvature tensor at $P_0$. 
\begin{equation*}
R_{\alpha\beta\gamma\delta}=\eps_{\delta}R_{\alpha\beta\gamma}^{\ \ \ \ \delta}=
\eps_{\delta}(\Gamma_{\alpha\gamma;\beta}^{\delta}-\Gamma_{\beta\gamma;\alpha}^{\delta}
)=\frac{1}{2}\big[g_{\beta\gamma;\alpha\delta}+g_{\alpha\delta;\beta\gamma}
-g_{\alpha\gamma;\beta\delta}-g_{\beta\delta;\alpha\gamma}\big].
\end{equation*}
We use this formula to compute the curvature components. 
\begin{thm} The only possibly non-vanishing components of the curvature tensor are listed below. 
\begin{enumerate}
\item $R(e_i,e_j,e_i,e_j)=-R(e_i,e_j,e_j,e_i)=4\eps_i\eps_j$;
\item  
$R(f_i,f_j,f_i,f_j)=-R(f_i,f_j,f_j,f_i)=4\eps_i\eps_j$. 
\item $R(e_i,e_j,f_k,f_l)=R(f_k,f_l,e_i,e_j)=-<x_i\overline{x_l},x_j\overline{x_k}>+<x_j\overline{x_l},x_i\overline{x_k}>$;
\item $R(e_i,f_j,e_k,f_l)=R(f_i,e_j,f_k,e_l)=<x_i\overline{x_j},x_k\overline{x_l}>$ 
and $R(e_i,f_j,f_l,e_k)=$ $R(f_i,e_j,e_l,f_k)=
-<x_i\overline{x_j},x_k\overline{x_l}>$.
\end{enumerate}
\end{thm}

The proof of this theorem is analogous to the proof of theorem \ref{BIGTHM} and will be omitted. Instead, we discuss a component-free description of the Riemann curvature tensor of $\oct' P^2$. For that purpose note that one could summarize parts (1) and (2) of the theorem by  
\begin{alignat*}{5}
&R(e_i,e_j,e_k,e_l)=R(f_i,f_j,f_k,f_l)\\
=&-4<x_i,x_l><x_j,x_k>+4<x_j,x_l><x_i,x_k>.
\end{alignat*}
Identify the tangent space at $P_0[1,0,0]$ with pairs of  para-octonions $(a,b)$ using the correspondence
$$(a,b)=(\textstyle{\Sigma} a_ix_i,\textstyle{\Sigma} b_ix_i)\leftrightarrow\textstyle{\Sigma}a_ie_i+\textstyle{\Sigma}b_if_i.$$
The following formula follows from the $\mathbb{R}$-multilinearity of the curvature tensor, the inner-product, the multiplication and the conjugation of para-octonions. 

\begin{cor}
\label{R'}
The Riemann curvature tensor of $\oct' P^2$ at the point $P_0[1,0,0]$ is given by
\begin{alignat*}{5}
R\big( (a,b),(c,d),(e,f),(g,h)\big)=&4<a,e><c,g>-4<c,e><a,g>\\+&4<b,f><d,h>-4<d,f><b,h>\\
-&<e\bar{d},g\bar{b}>+<e\bar{b},g\bar{d}>-<c\bar{f}, a\bar{h}>+<a\bar{f},c\bar{h}>\\
-&<a\bar{d}-c\bar{b},g\bar{f}-e\bar{h}>.
\end{alignat*}
\end{cor}

\section{Indefinite octonionic projective plane $\oct P^{(1,1)}$}\label{(1,1)}\label{op11section}

To define the indefinite octonionic projective plane $\oct P^{(1,1)}$ we again
consider the relation $\sim$ on $\oct^3$ 
defined by 
$$[a,b,c]\sim [d,e,f] \Longleftrightarrow a=d\lambda, b=e\lambda, c=f\lambda 
\text{\ \ for some\ \ } \lambda \in \oct-\{0\}.$$ 
This time we start with the sets
\begin{alignat*}{5}
U_1:=&\{(1,u,v)\in\mathbb{O}^3\big{|} 
1+|u|^2-|v|^2>0\},\\ 
U_2:=&\{(u,1,v)\in\mathbb{O}^3\big{|} |u|^2+1-|v|^2>0\}
\end{alignat*}
and study the restriction of $\sim$ to the union  $\mathcal{U}:=U_1\cup U_2$. \begin{defn}
The indefinite octonionic projective plane is the set of equivalence classes of $\mathcal{U}$ by the equivalence relation $\sim$.
$$\oct P^{(1,1)}=\mathcal{U}/_{\sim}$$
\end{defn}
We can topologize $\oct P^{(1,1)}$ using the sets $U_i$; with this topology $\oct P^{(1,1)}$ becomes a 16-dimensional manifold (see also theorem \ref{mfld}). In fact, we can use the Van Kampen theorem to show that  $\oct P^{(1,1)}$ is  a simply connected manifold. The sets 
$$U_1/_{\sim}\approx U_2/_{\sim}\approx \big\{(u,v)\in \Real^{16}\big| |u|^2-|v|^2>-1\big\}$$  are simply connected due to the homotopy $H\big((u,v),t\big)=(u,tv)$; this homotopy proves that $U_i/_{\sim}\simeq \Real^8$.
The same homotopy shows that 
$$U_1/_{\sim}\cap U_2/_{\sim}=\big\{[1,u,v]\big|1+|u|^2-|v|^2>0,\ u\ne 0\big\}\simeq \oct-\{0\}\simeq S^7.$$
Hence $U_1/_{\sim}\cap U_2/_{\sim}$ is simply connected. 
Since the sets $U_1/_{\sim}, U_2/_{\sim}$ and $U_1/_{\sim}\cap U_2/_{\sim}$ are all simply connected so is $\oct P^{(1,1)}=U_1/_{\sim}\cup U_2/_{\sim}$.

We now introduce the metric structure on $\oct P^{(1,1)}$. Consider 
\begin{equation}\label{7.?}
ds^2=\frac{|du|^2(1-|v|^2)-|dv|^2(1+|u|^2)
+2Re[(u\bar{v})(dv d\bar{u})]}{(1+|u|^2-|v|^2)^2}\ 
\end{equation}
on the charts $U_1/_{\sim}$ and $U_2/_{\sim}$.
The compatibility of the two metrics  on $U_1/_{\sim}\cap U_2/_{\sim}$ can be verified using methods analogous to those of theorem \ref{thm0}. 

To establish non-degeneracy of our metric we again consider the metric tensor components. As in section \ref{para}, we use the standard orthonormal basis $\{x_1,...,x_8\}$ of $\oct$, the coordinate frame
$\{e_1,...,e_8,f_1,...,f_8\}$ where 
$$e_i:=\partial_i,\ f_i:=\partial_{i+8},\ 1\le i\le 8$$ 
and the $8\times 8$ matrix $A$ with components $A_{ij}=-<(u\bar{v})x_j,x_i>$.

Using this notation we can write our metric tensor $g$
as
$$\mathcal{M}=
\frac{1}{(1+|u|^2-|v|^2)^2}
\left[\begin{array}{cc}(1-|v|^2)Id & -A \\ -A^T & -(1+|u|^2)Id\end{array}\right].$$
An argument analogous to lemma \ref{kinda_orth} shows that $AA^T=A^TA=|u|^2|v|^2Id$. 

\begin{prop}\label{signature2}
The expression (\ref{7.?})  defines a non-degenerate metric on $\oct P^{(1,1)}$ of signature (8,8).
\end{prop}
\begin{proof}
To show that $\mathcal{M}$ is non-degenerate suppose there exists 
$\vec{v}=(\vec{r},\vec{s})$ such that  $\mathcal{M}\vec{v}=0$. 
The vectors $\vec{r}$, $\vec{s}$ then must satisfy
$$\begin{cases}
(1-|v|^2)\vec{r}-A\vec{s}=0 & \cr
-A^T\vec{r}-(1+|u|^2)\vec{s}=0 & 
\end{cases}
$$
which, using $A^TA=AA^T=|u|^2|v|^2Id$,  reduces to
$$(1+|u|^2-|v|^2)\vec{r}=0 \text{\ \ and\ \ }
(1+|u|^2-|v|^2)\vec{s}=0.$$
Since by assumption $1+|u|^2-|v|^2>0$, we see that $\vec{r}=\vec{s}=0$ i.e. that the inner product defined by $\mathcal{M}$  is non-degenerate. 

The signature of a non-degenerate metric is the same at every point. Thus we only need to consider the signature at one point. For convenience let us consider the point $[1,0,0]$; at this point $ds^2=|du|^2-|dv|^2$. Therefore our metric is of signature (8,8). 
\end{proof}

For symmetry reasons it would be good to know the form of our metric on 
$$U_3/_{\sim}=\big\{[1,u,v]\in U_1/_{\sim}\big| v\ne 0\big\}\cup \big\{[u,1,v]\in U_2/_{\sim}\big|v\ne 0\big\};$$
this set can be coordinatized using 
$$U_3/_{\sim}=\big\{[x,y,1]\big|\ |x|^2+|y|^2-1>0\big\}.$$
Due to the rational nature of the metric and the transition functions, it is enough to express the metric on $U_1/_{\sim}\cap U_3/_{\sim}$ using the coordinates arising from $U_3/_{\sim}$. In other words, in order to find the metric in terms of the coordinates on $U_3/_{\sim}$ we need to set $u=yx^{-1}$ and $v=x^{-1}$ in  (\ref{7.?}). A computation which is very analogous to the one we performed in the proof of theorem 
\ref{thm0} yields
\begin{equation}\label{7.???}
ds^2=\frac{|dx|^2(|y|^2-1)+|dy|^2(|x|^2-1)
-2Re[(x\bar{y})(dy d\bar{x})]}{(|x|^2+|y|^2-1)^2}.\ 
\end{equation}

We now take a look at some of the isometries of $\oct P^{(1,1)}$. 
For $r^2-|\lll |^2=1$ consider the ``indefinite reflections" 
$\rho_{r,\lll}$ on $U_1/_{\sim}$ (or $U_2/_{\sim}$) defined by
\begin{equation*}
u'= -ru+\lll v\text{\ \ and\ \ } v'= -\bar{\lll}u+rv.
\end{equation*}
The images of these maps are contained in $U_1/_{\sim}$ since $|u'|^2-|v'|^2=|u|^2-|v|^2$. 
A lengthy computation similar to the one of proposition \ref{isometry1} shows that the maps $\rho_{r,\lll}$ are (local) isometries satisfying $\rho_{r,\lll}^2=Id$. 
As expected, these ``reflections" extend to globally defined maps on $\oct P^{(1,1)}$. The {\it formal} extension of $\rho_{r,\lll}$ to $U_2/_{\sim}$ can be computed using methods of  proposition \ref{extensions}:
\begin{equation}\label{rationalagain}
[x,1,z]\mapsto 
 \Big[\ \frac{-rx+(x\bar{z})\bar{\lll}}
 {|r-\lll z|^2},\ 1,\ \frac{r\bar{\lll}-\bar{\lll}\bar{z}\bar{\lll}-r^2z+r|z|^2\bar{\lll}}{|r-\lll z|^2}\ \Big].
 \end{equation}
The last expression may not be defined on all of $U_2/_{\sim}$. For this reason we consider the set
\begin{equation*}
U_2':=\Big\{[x,1,z]\in U_2/_{\sim}\big| r-\lll z\ne 0\Big\}.
\end{equation*} 
Note that  if $r-\lll z=0$ for some $[x,1,z]\in U_2/_{\sim}$ 
then $\lll\ne 0$, 
$$0<|x|^2+1-|z|^2=|x|^2+1-\frac{r^2}{|\lll|^2}=
|x|^2-\frac{1}{|\lll|^2}$$
and consequently $x\ne 0$. Therefore, we have an open cover 
\begin{equation}
\oct P^{(1,1)}=U_1/_{\sim}\cup U_2',
\end{equation}
and 
the methods of proposition \ref{extensions} apply to our current situation. Hence  the maps $\rho_{r,\lll}$ extend to global isometries of $\oct P^{(1,1)}$.

There is another kind of isometry on $\oct P^{(1,1)}$. It arises from ``Euclidean reflections" $\mathcal{R}_{r,\lll}$ on $U_3/_{\sim}$:
\begin{equation}\label{7.6}
[x,y,1]\mapsto [x',y',1], \text{\ \ where\ \ }x'= rx+\lll y,\ y'= \bar{\lll}x-ry, \text{\ \ and\ \ }  r^2+|\lll|^2=1.
\end{equation}
The map  $\mathcal{R}_{r,\lll}$ is well-defined due to $|x'|^2+|y'|^2=|x|^2+|y|^2$. This identity also shows that $\mathcal{R}_{r,\lll}$ is an  isometry if and only if it preserves 
$$|dx|^2|y|^2+|dy|^2|x|^2-2Re[(x\bar{y})(dyd\bar{x})].$$
The condition was already checked in the proof of the proposition \ref{isometry1} and so  $\mathcal{R}_{r,\lll}$ indeed is an isometry. We can also easily see that this map is an involution.

The same computation as in proposition \ref{extensions} shows that {\it formal} extension of $\mathcal{R}_{r,\lll}$ to $U_1/_{\sim}$ is 
\begin{equation}\label{ext.(1,1)again}
[1,u,v]\mapsto 
\begin{cases}
       \Big[1,
       \ \frac{r\bar{\lll}-r^2u+\bar{\lll}\bar{u}\bar{\lll}-r|u|^2\bar{\lll}}{|r+\lll u|^2},\ \frac{rv+(v\bar{u})\bar{\lll}}{|r+\lll u|^2}\ \Big] 
       &  \text{\  for\  }[1,u,v]\in U_1', \cr
        &       \cr
         \Big[\ \frac{r\lll+\lll u\lll-r^2\bar{u}-r|u|^2\lll}{|\bar{\lll}-ru|^2},\ \frac{v\lll-r(v\bar{u})}{|\bar{\lll}-ru|^2},\ 1\Big] & \text{\  for\  }[1,u,v]\in U_1'' ,
 \end{cases}
\end{equation}
where $U_1'=\{[1,u,v]\in U_1/_{\sim}\big|\ r+\lll u\ne 0\}$ and 
$U_1''=\{[1,u,v]\in U_1/_{\sim}\big|\ \bar{\lll}-ru\ne 0\}$. Since $r+\lll u =0$ and  $\bar{\lll}-ru = 0$ implies 
$r^2+|\lll|^2=r(r+\lll u)+\lll (\bar{\lll}-ru)=0$, we see that $\{U_1', U_1''\}$ is an open cover of $U_1/_{\sim}$. The map $\mathcal{R}_{r,\lll}$ formally extends to $U_2/_{\sim}$ in a similar fashion and we omit the details.  

We now explain why the map $\mathcal{R}_{r,\lll}$ is well-defined. 
Consider a point $$[1,u,v]=[x,y,1]\in U_1'\cap U_3/_{\sim}.$$ Its image $[1,u',v']=[x',y',1]$ satisfies 
\begin{alignat*}{5}
1+|u'|^2-|v'|^2=& |v'|^2 (|x'|^2+|y'|^2-1)=
|v'|^2(|x|^2+|y|^2-1)\\
= & \frac{|v'|^2}{|v|^2}(1+|u|^2-|v|^2)
= \frac{1}{|r+\lll v|^2}(1+|u|^2-|v|^2)
\end{alignat*}
and thus the formal image of $[1,u,v]\in U_1'$ indeed is an element of $U_1/_{\sim}$ {\it for all} $[1,u,v]\in U_1'$.  A similar argument applies in all the remaining cases ($U_1'', U_2', U_2''$).
This globally defined map is an involution and an isometry due to the rational nature of the map (\ref{7.6}), the rational nature of our transition functions, and the fact that  
$\mathcal{R}_{r,\lll}$ is an involutive isometry on $U_3/_{\sim}$.

\begin{thm} $\oct P^{(1,1)}$ is homogeneous. 
\end{thm}
\begin{proof}
It suffices to show that the point $[1,0,0]$ can be taken to any of the points $[1,a,b]$ with $1+|a|^2-|b|^2>0$. In fact, without loss of generality we may assume that 
$|b|<1$. To see this consider the natural extension of the ``Euclidean reflection" $\mathcal{R}_{r,\lll}$ (see (\ref{ext.(1,1)again})) with $r=\frac{1}{\sqrt{1+|a|^2}}$, $\lll=\frac{\bar{a}}{\sqrt{1+|a|^2}}$. This map takes $[1,a,b]$ to $[1,a',b']$ where
$$b'=\frac{\frac{b}{\sqrt{1+|a|^2}}+\frac{(b\bar{a})a}{\sqrt{1+|a|^2}}}{\big|\frac{1}{\sqrt{1+|a|^2}}+\frac{\bar{a}a}{\sqrt{1+|a|^2}}\big|^2}=\frac{b\sqrt{1+|a|^2}}{1+|a|^2}=\frac{b}{\sqrt{1+|a|^2}}.$$
As $1+|a|^2-|b|^2>0$ we see that $\frac{|b|^2}{1+|a|^2}<1$, i.e $|b'|<1$.

To build an isometry taking $[1,0,0]$ to $[1,a,b]$ with $|b|<1$ first consider a real number $t_0$ such that 
$$\tan t_0= \frac{-|a|}{\sqrt{1-|b|^2}}.$$
Let $a_0$ be the unit octonion $a_0=\frac{\bar{a}}{|a|}$ (if $a=0$ just consider {\it any} unit octonion $a_0$). 
The ``Euclidean reflection" $\mathcal{R}_{\cos t_0, \sin t_0 a_0}$ maps 
$$[1,0,0]\mapsto [1,\frac{-a}{\sqrt{1-|b|^2}},0].$$
We now consider the ``indefinite reflection" on $U_2/_{\sim}$:
\begin{equation*}
[x,1,z]\mapsto [-rx+\lll z, 1, -\bar{\lll}x+rz],\ \text{\ with\ } r^2-|\lll|^2=1.
\end{equation*}
The natural extension of this map to (a  certain subset of) $U_1/_{\sim}$ is given by 
\begin{equation}\label{indef.refl}
[1,u,v]\mapsto 
 \Big[1, \ \frac{-ru+(u\bar{v})\bar{\lll}}
 {|r-\lll v|^2},\ \frac{r\bar{\lll}-\bar{\lll}\bar{v}\bar{\lll}-r^2v+r|v|^2\bar{\lll}}{|r-\lll v|^2}\ \Big],
\end{equation}
see (\ref{rationalagain}) for details.
Choose $r=\frac{1}{\sqrt{1-|b|^2}}$ and $\lll=\frac{\bar{b}}{\sqrt{1-|b|^2}}$. 
The map (\ref{indef.refl}) takes 
$[1,\frac{-a}{\sqrt{1-|b|^2}},0]$ to 
$$\Big[1,\ -\frac{\ \frac{-a}{\sqrt{1-|b|^2}}}{r},\  \frac{\bar{\lll}}{r}\Big]=[1,a,b].$$
Therefore, $\oct P^{(1,1)}$ is a homogeneous manifold. 
\end{proof}

The curvature of $\oct P^{(1,1)}$ can also be studied using methods of section \ref{curvsect}. As $\oct P^{(1,1)}$ is homogeneous we may restrict our attention to the point $P_0[1,0,0]$.

Our metric at $P_0$ takes the form of   $$g|_{P_0}=\text{\ diag\ } (1,1,1,1,-1,-1,-1,-1);$$ 
for notational simplicity  we shall use 
$\eps_i=\begin{cases}
1 &\text{\ if\ } i\le 4,\cr
-1 &\text{\ if\ }i\ge 5
\end{cases}$.
To study the jets of the metric at $P_0$ we first find 
the appropriate quadratic approximations using Maclaurin series. Since $$\frac{1}{(1+|u|^2-|v|^2)^2}=1-2(|u|^2-|v|^2)+O\Big(\big(|u|^2-|v|^2\big)^2\Big), $$
we have:
\begin{alignat*}{5} 
-\frac{1+|u|^2}{(1+|u|^2-|v|^2)^2}&\approx
-1+2(|u|^2-|v|^2)-|u|^2=-1+|u|^2-2|v|^2\\ 
\frac{1-|v|^2}{(1+|u|^2+|v|^2)^2}&\approx1-2(|u|^2-|v|^2)-|v|^2=1-2|u|^2+|v|^2\\
\frac{<(u\bar{v})x_j,x_i>}{(1+|u|^2-|v|^2)^2}&\approx  <(u\bar{v})x_j,x_i>.
\end{alignat*}
We can now easily see that the first jets of the metric vanish at $P_0$, while the only 
possibly non-vanishing second jets are: 
\begin{itemize}
\item $e_je_jg(e_i,e_i)=f_jf_jg(f_i,f_i)=-4$;
\item $f_jf_jg(e_i,e_i)=e_je_jg(f_i,f_i)=2$;
\item $e_lf_kg(e_i,f_j)=<x_l\overline{x_k}, x_i\overline{x_j}>$.
\end{itemize}
Therefore, at $P_0$ we have 
\begin{equation*}
R_{\alpha\beta\gamma\delta}=\eps_{\delta}R_{\alpha\beta\gamma}^{\ \ \ \ \delta}=
\eps_{\delta}(\Gamma_{\alpha\gamma;\beta}^{\delta}-\Gamma_{\beta\gamma;\alpha}^{\delta}
)=\frac{1}{2}\big[g_{\beta\gamma;\alpha\delta}+g_{\alpha\delta;\beta\gamma}
-g_{\alpha\gamma;\beta\delta}-g_{\beta\delta;\alpha\gamma}\big]
\end{equation*}
which gives us the following theorem.

\begin{thm} \label{BIGTHM(1,1)} The only possibly non-vanishing components of the curvature tensor are listed below. 
\begin{enumerate}
\item $R(e_i,e_j,e_i,e_j)=-R(e_i,e_j,e_j,e_i)=4$;
\item $R(f_i,f_j,f_i,f_j)=-R(f_i,f_j,f_j,f_i)=4$;
\item $R(e_i,e_j,f_k,f_l)=R(f_k,f_l,e_i,e_j)=<x_i\overline{x_l},x_j\overline{x_k}>-<x_j\overline{x_l},x_i\overline{x_k}>$;
\item $R(e_i,f_j,e_k,f_l)=R(f_i,e_j,f_k,e_l)=-<x_i\overline{x_j},x_k\overline{x_l}>$ and $R(e_i,f_j,f_l,e_k)=R(f_i,e_j,e_l,f_k)=
<x_i\overline{x_j},x_k\overline{x_l}>$.\end{enumerate}
\end{thm}

It should be noted that parts (1) and (2) of theorem \ref{BIGTHM(1,1)} can be summarized as 
\begin{alignat*}{5}
&R(e_i,e_j,e_k,e_l)=R(f_i,f_j,f_k,f_l)\\
=&-4<x_i,x_l><x_j,x_k>+4<x_j,x_l><x_i,x_k>.
\end{alignat*}
Using $\mathbb{R}$-multilinearity of the expressions in theorem \ref{BIGTHM(1,1)} we obtain the following component-free description of the Riemann curvature tensor. \begin{cor} The Riemann curvature tensor of $\oct P^{(1,1)}$ at the point $P_0[1,0,0]$ is given by
\begin{alignat*}{5}
R\big( (a,b),(c,d),(e,f),(g,h)\big)=&4<a,e><c,g>-4<c,e><a,g>\\+&4<b,f><d,h>-4<d,f><b,h>\\
+&<e\bar{d},g\bar{b}>-<e\bar{b},g\bar{d}>+<c\bar{f}, a\bar{h}>-<a\bar{f},c\bar{h}>\\
+&<a\bar{d}-c\bar{b},g\bar{f}-e\bar{h}>,
\end{alignat*}
where we  identified the tangent space at $P_0[1,0,0]$ with pairs of octonions $(a,b)$ according to  
$$(a,b)=(\textstyle{\Sigma} a_ix_i,\textstyle{\Sigma} b_ix_i)\leftrightarrow\textstyle{\Sigma}a_ie_i+\textstyle{\Sigma}b_if_i.\qed$$
\end{cor}

\section{The octonionic hyperbolic plane $\oct H^2$}

In this section we study our final example: 
the octonionic hyperbolic plane, i.e. 
the hyperbolic dual of the octonionic projective plane.
\begin{defn}\label{defoh2}
The octonionic hyperbolic plane is the set 
$$\oct H^2=\big\{(u,v)\in\oct^2\big| |u|^2+|v|^2<1\big\}$$ equipped with the metric 
\begin{equation}\label{hypmetric}
ds^2=\frac{|du|^2(1-|v|^2)+|dv|^2(1-|u|^2)
+2Re[(u\bar{v})(dv d\bar{u})]}{(1-|u|^2-|v|^2)^2}.
\end{equation}
\end{defn}

Note that we can define $\oct H^2$ equivalently as $\big\{[1,u,v]\in \oct^3\big| 1-|u|^2-|v|^2>0\big\}$. This  interpretation makes the connection between $\oct H^2$ and the previous examples more clear. 

The metric in (\ref{hypmetric}) is positive definite. To see this we study the metric components with respect to frame $\{e_1, \ldots, e_8, f_1,\ldots, f_8\}$ where $e_i:=\partial_i$, $f_i:=\partial_{i+8}$. The metric tensor of (\ref{hypmetric}) has the following matrix representation with respect to $\{e_1,\ldots,f_8\}$:
$$\mathcal{M}=
\frac{1}{(1-|u|^2-|v|^2)^2}
\left[\begin{array}{cc}(1-|v|^2)Id & A \\ A^T & (1-|u|^2)Id\end{array}\right].$$
Here $A_{ij}=<(u\bar{v})x_j,x_i>$, where $\{x_1,\ldots,x_8\}$ is the standard orthonormal basis for $\oct$. 
Using the methods of section \ref{para} we see that the non-degenerateness of (\ref{hypmetric}) reduces to showing that the system 
\begin{equation*}
\begin{cases}
(1-|v|^2)\vec{r}+A\vec{s}=0 & \cr
A^T\vec{r}+(1-|u|^2)\vec{s}=0 & 
\end{cases}
\end{equation*}
has only trivial solutions. As $AA^T=A^TA=|u|^2|v|^2Id$ this system yields
$$\Big(-(1-|v|^2)(1-|u|^2)+|u|^2|v|^2\Big)\vec{s}=0 
\text{\ \ and\ \ }
\Big(|u|^2|v|^2-(1-|v|^2)(1-|u|^2)\Big)\vec{r}=0.$$
Since $|u|^2|v|^2-(1-|v|^2)(1-|u|^2)=-1+|u|^2+|v|^2<0$ we see that $\vec{r}=\vec{s}=0$, i.e. that the metric (\ref{hypmetric}) is non-degenerate. To establish it is in fact positive definite we now only need to check it is positive definite at one particular point. For convenience we may consider $P_0$ with coordinates $(0,0)$ where our metric takes the form of $|du|^2+|dv|^2$. As $|du|^2+|dv|^2$ is positive definite, so is the metric (\ref{hypmetric}). 
We have therefore proven the following:
\begin{thm}
$\oct H^2$ of definition \ref{defoh2} is a $16$-dimensional simply connected Riemannian manifold. 
\end{thm}

For symmetry reasons it would be good to know the form of our metric on 
$$U_2/_{\sim}=\big\{[1,u,v]\in \oct H^2 \big| u\ne 0\big\} \text{\ \ \ and\ \ \ }
U_3/_{\sim}=\big\{[1,u,v]\in \oct H^2 \big| v\ne 0\big\};$$
these sets can be coordinatized using 
$$U_2/_{\sim}=\big\{[x,1,z]\in \big| |x|^2-1-|z|^2> 0\big\} \text{\ \ \ and\ \ \ }
U_3/_{\sim}=\big\{[x,y,1]\big|\ |x|^2-|y|^2-1>0\big\}.$$
Setting $u=x^{-1}$ and $v=zx^{-1}$ into (\ref{hypmetric}) gives us the metric on $U_2/_{\sim}$. We obtain:
\begin{equation}\label{hypmetric-1}
ds^2=\frac{|dx|^2(1+|z|^2)+|dz|^2(|x|^2-1)
-2Re[(x\bar{z})(dz d\bar{x})]}{(|x|^2-1-|z|^2)^2}.\ 
\end{equation}
Similarly, by replacing  $u=yx^{-1}$ and $v=x^{-1}$ in  (\ref{hypmetric}) we see that the metric on $U_3/_{\sim}$ is given by 
\begin{equation}\label{hypmetric-2}
ds^2=\frac{|dx|^2(1+|y|^2)+|dy|^2(|x|^2-1)
-2Re[(x\bar{y})(dy d\bar{x})]}{(|x|^2-|y|^2-1)^2}.\ 
\end{equation}
All the computations involved are completely analogous to the one in the proof of theorem \ref{thm0}.

There are at least two kinds of isometries on $\oct H^2$: those arising from ``Euclidean reflections" on $\oct H^2$ and those arising from ``indefinite reflections" on $U_2/_{\sim}$ and $U_3/_{\sim}$.
\begin{itemize}
\item ``Euclidean reflections" $\mathcal{R}_{r,\lll}$ on $\oct H^2$ take the form $[1, u, v]\mapsto [1, u', v']$, where $u'= ru+\lll v,\ v'= \bar{\lll}u-rv$ and $r^2+|\lll|^2=1$. These are well-defined on $\oct H^2$ since $|u'|^2+|v'|^2=|u|^2+|v|^2<1$. It is easy to verify that $\mathcal{R}_{r,\lll}^2=Id$. The maps $\mathcal{R}_{r,\lll}$ are isometries 
since $|du'|^2+|dv'|^2=|du|^2+|dv|^2$ and 
since 
\begin{alignat*}{5}
&-|du'|^2|v'|^2-|dv'|^2|u'|^2+2Re[(u'\bar{v'})(dv'd\bar{u'})]\\
=&  
-|du|^2|v|^2-|dv|^2|u|^2+2Re[(u\bar{v})(dvd\bar{u})]
\end{alignat*} due to the proof of proposition \ref{isometry1}.
\item The ``indefinite reflections" $\rho_{r,\lll}$ on $U_2/_{\sim}$ take the form 
\begin{equation*}
x'= -rx+\lll z,\ z'= -\bar{\lll}x+rz, \text{\ \ \ where\ \ \ }r^2-|\lll |^2=1.
\end{equation*}
As usual, we have $\rho_{r,\lll}^2=Id$ along with $|x'|^2-|z'|^2=|x|^2-|z|^2$. The latter ensures that $\rho_{r,\lll}:U_2/_{\sim}\to U_2/_{\sim}$ is well-defined.  The maps $\rho_{r,\lll}$ are isometries 
since $|dx'|^2-|dz'|^2=|dx|^2-|dz|^2$ and 
since 
\begin{alignat*}{5}
&|dx'|^2|z'|^2+|dz'|^2|x'|^2-2Re[(x'\bar{z'})(dz'd\bar{x'})]\\
=&  
|dx|^2|z|^2+|dz|^2|x|^2-2Re[(x\bar{z})(dzd\bar{x})].
\end{alignat*} 
This last identity has not been proven in this paper per se, but it is the essential part of the statement that the maps $\rho_{r,\lll}$ of section \ref{(1,1)} are isometries. One could define analogous isometries on $U_3/_{\sim}$.
\end{itemize}

We now discuss the extension of $\rho_{r,\lll}$ to $\oct H^2$. As in equation (\ref{rationalagain}), the {\it formal} extension of $\rho_{r,\lll}$ is 
\begin{equation*}
(u,v)\mapsto (u',v') \text{\ \ with\ \ }
u'=\frac{-ru+(u\bar{v})\bar{\lll}}
 {|r-\lll v|^2},\ 
v'=\frac{r\bar{\lll}-\bar{\lll}\bar{v}\bar{\lll}-r^2v+r|v|^2\bar{\lll}}{|r-\lll v|^2}.
 \end{equation*}
Note that the expressions for $u'$ and $v'$ are well-defined since $-r+\lll v\ne 0$
for all $(u,v)\in \oct H^2$. This is because if $-r+\lll v=0$ then $\lll\ne 0$, $v=r\lll^{-1}$ and $$|v|^2=\frac{r^2}{|\lll|^2}=\frac{1+|\lll|^2}{|\lll|^2}>1.$$
To see that $|u'|^2+|v'|^2<1$ we first study the points $[1,u,v]=[x,1,z]\in U_2/_{\sim}$. We have:
\begin{alignat*}{5}
1-|u'|^2-|v'|^2&=|v'|^2\big(|x'|^2-1-|z'|^2\big)=|v'|^2\big(|x|^2-1-|z|^2\big)\\
&=
\frac{|v'|^2}{|v|^2}\big(1-|u|^2-|v|^2\big)=\frac{1}{|r-\lll v|^2}\big(1-|u|^2-|v|^2\big).
\end{alignat*}
Due to its rational nature this equality holds for all $(u,v)\in \oct H^2$ and therefore $|u'|^2+|v'|^2<1$. Now that we know that the extensions of $\rho_{r,\lll}$ are well-defined we may use the standard argument involving the rational nature of $\rho_{r,\lll}$ to show that these maps are isometries of $\oct H^2$.

\begin{prop}
$\oct H^2$ is homogeneous. 
\end{prop}

\begin{proof}
It is enough to show that the point $P_0$ with coordinates $(0,0)$ (i.e $[1,0,0]$) can be taken to any point $(a,b)$, $0<|a|^2+|b|^2<1$ via an isometry.  

Let $R:=\sqrt{|a|^2+|b|^2}$. We first show that $P_0$ can be taken to $(0,R)$ via an ``indefinite reflection" $\rho_{r,\lll}$ arising from $U_2/_{\sim}$. Let $t\in \mathbb{R}$ be such that $R=\tanh t$; such a $t$ exists as $R<1$. The isometry $\rho_{r,\lll}$ with  
$r=\cosh t$ and $\lll = \sinh t $ takes 
$P_0\mapsto (0, \tanh t)=(0,R).$

We continue by considering two cases.

{\it Case 1,} when $a=0$. Then $b\ne 0$ and we may consider the composition of ``Euclidean reflections" $\mathcal{R}_{0,\frac{\bar{b}}{|b|}}\circ \mathcal{R}_{0,1}$. We have
$$(0,R)=(0,|b|)\mapsto (|b|,0)\mapsto (0,b)=(a,b).$$

{\it Case 2,} when $a\ne 0$. We consider the following composition of ``Euclidean reflections":
$$\mathcal{R}_{0,1}\circ\mathcal{R}_{\frac{|a|}{R}, -\frac{b\bar{a}}{R|a|}}\circ\mathcal{R}_{0,-\frac{\bar{a}}{|a|}}\circ \mathcal{R}_{0,1}.$$
The effect of this composition on $[1,0,R]$ is easily seen to be
$$(0,R)\mapsto(R,0)\mapsto\Big(0,-\frac{Ra}{|a|}\Big)\mapsto(b,a)\mapsto (a,b).$$
This completes the proof that $\oct H^2$ is homogenous. 
\end{proof}

The point $P_0$ with coordinates $(0,0)$ is also convenient for computing the curvature of $\oct H^2$. The metric at $P_0$ is Euclidean, while the jets of the metric at $P_0$ can be computed from the following quadratic approximations:
\begin{alignat*}{5} 
\frac{1-|u|^2}{(1-|u|^2-|v|^2)^2}&\approx1+2(|u|^2+|v|^2)-|u|^2=1+|u|^2+2|v|^2\\ 
\frac{1-|v|^2}{(1-|u|^2-|v|^2)^2}&\approx1+2(|u|^2+|v|^2)-|v|^2=1+2|u|^2+|v|^2\\
\frac{<(u\bar{v})x_j,x_i>}{(1-|u|^2-|v|^2)^2}&\approx  <(u\bar{v})x_j,x_i>.
\end{alignat*}
Compare these with the approximations (\ref{approx1})-(\ref{approx3}); it is easy to see that the first jets of (\ref{hypmetric}) at $P_0$ vanish and that the second jets at $P_0$ are the exact negatives of those from lemma \ref{jets}. Hence the following result.
\begin{thm}\label{oh2}
The Riemann curvature tensor of $\oct H^2$ at the point $P_0(0,0)$ is given by the negative of the curvature tensor of $\oct P^2$ at $[1,0,0]$ (see corollary  \ref{R}).\qed
\end{thm}

\section{Identification with classical models}
In this section we use the classification results of E. Garcia-Rio, D. N. Kupeli and R. Vazques-Lorenzo \cite{spanishguys} regarding {\it semi-Riemannian special Osserman manifolds} to identify our (para-)octonionic projective planes with the projective planes defined using  exceptional Lie groups \cite{Lie Algebras, wolf}. We start by explaining the context of special Osserman manifolds in more detail. 

By the {\it Jacobi operator} (at a point $P$) of a semi-Riemannian manifold $M$ we mean the family of self-adjoint operators (on the tangent space at $P$) defined by the Riemann curvature tensor $R$ as follows:
$$\J_v(x)=R(v,x)v,\ \ x,v\in T_PM.$$
In the case of a {\it (locally) isotropic} manifold, i.e. a manifold $M$ such that for any $P\in M$ and any two non-zero tangent vectors $v,w$ at $P$ with $g(v,v)=g(w,w)$ there exists a (local) isometry preserving $P$ whose differential takes $v$ to $w$, the spectrum of the operator $\J_v$ is independent of the choice of unit spacelike (resp. unit timelike, non-zero null) vector $v$ at $P$.  The converse is not true in general, as evidenced by the para-complex projective plane (see \cite{Lie Algebras, elbook, wolf}). However, in Riemannian geometry  Osserman \cite{osserman} conjectured that the converse holds. This has been proven for manifolds of dimension other than 16 by Nikolayevski (see \cite{nikolayevski}).

Osserman's conjecture initiated the study of the so-called {\it Osserman manifolds}, i.e. manifolds for which the spectrum of the Jacobi operator is constant over the (spacelike or timelike) unit sphere bundles. Classification of Osserman manifolds is a hard problem and smaller classes of Osserman manifolds are considered instead. In \cite{spanishguys} the authors study {\it special Osserman manifolds} which are characterized by the following conditions.

Let $v\in T_PM$ be unit. Since $\J_v(v)=0$ the only interesting part of the spectrum comes from the restriction $\J_v:v^\perp\to v^\perp$. Note that for Osserman manifolds the spectrum of $J_v$ changes sign depending on whether $v$ is spacelike or timelike.
\begin{itemize}
\item {\it (Condition I)} The operator $\J_v:v^\perp\to v^\perp$ is diagonalizable with exactly $2$ non-zero eigenvalues $\eps_v\lll$ and $\eps_v\mu$; here we use $\eps_v=g(v,v)$ to account for the sign difference in the spectrum.
\end{itemize}
The remaining conditions concern the space $E_\lll(v):=\mathrm{span}\{v\}\oplus\ker\{\J_v-\eps_v\lll Id\}$ and 
the $\mu$ eigenspace of the Jacobi operator. \begin{itemize}
\item {\it (Condition II)} If $v,w$ are unit and $w\in E_\lll(v)$ then $E_\lll(v)=E_\lll(w)$;
\item {\it (Condition III)} If $v$ is a unit vector and $w\in\ker\{\J_v-\eps_v\mu Id\}$, then also $v\in\ker\{\J_w-\eps_w\mu Id\}$.
\end{itemize}

The classification result for special Osserman manifolds states the following.
\begin{thm}
The only complete and simply connected semi-Riemannian special Osserman manifolds are:
\begin{enumerate}
\item complex space forms with definite or indefinite metric tensor and para-complex space forms;
\item quaternionic space forms  with definite or indefinite metric tensor and 
para-quaternionic space forms;
\item octonionic projective plane with definite of indefinite metric tensor and the para-octonionic projective plane.
\end{enumerate}
\end{thm}

The basic examples of (para-)complex and (para-)quaternionic space forms are the (para-)complex and the (para-)quaternionic projective spaces (see \cite{elbook}). The hyperbolic duals of the Riemannian projective spaces are included in this classification as they arise from the negative definite projective spaces after the metric sign change. For example,  if $P(\mathbb{C}^{(n,1)})$ stands for the complex projective space of all spacelike lines in $\mathbb{C}^{(n,1)}$ the resulting 
Fubini-Study metric, 
$g_{FS}$, is negative definite. The hyperbolic dual $\mathbb{C}H^n$ of $\mathbb{C}P^n$ can be seen as 
$\big(P(\mathbb{C}^{(n,1)}), -g_{FS}\big)$.

The three categories of special Osserman manifolds can be distinguished by the multiplicities of the non-zero eigenvalues of the Jacobi operator. The Jacobi operator of a (para-)complex (resp. (para-)quaternionic) space form has one non-trivial eigenvalue of multiplicity $1$ (resp. $3$). In the case of (para-)octonionic planes this multiplicity is $7$. 

Our octonionic projective planes are simply connected  and homogenous, which in particular means they are complete. Therefore, to identify our projective planes with the standard models we only need to show they are special Osserman manifolds whose Jacobi operators have a non-trivial eigenvalue of multiplicity $7$. To see that $\oct' P^2$ indeed corresponds to the classical para-octonionic projective plane (and not to our other example with indefinite metric, $\oct P^{(1,1)}$) we prove that $\oct' P^2$ is {\it not} locally isotropic. This is sufficient due to Wolf's classification of locally isotropic semi-Riemannian manifolds \cite{wolf}.

\begin{thm}
$\oct P^2$, $\oct'P^2$ and $\oct H^2$ are special Osserman manifolds. 
\end{thm}

\begin{proof} It follows from theorem \ref{oh2} that it is enough to consider $\oct P^2$ and $\oct'P^2$. As these two manifolds are homogenous we may restrict our attention to the curvature tensors at $[1,0,0]$ where their metrics take the form 
$$g\big((a,b),(c,d)\big)=<a,c>+<b,d>.$$

Let $(c',d'):=\jac(c,d)$, where $\jac$ is the Jacobi operator at $[1,0,0]$ corresponding to a unit tangent vector $(a,b)$. 
We can express $c',d'$  using corollary \ref{R} and the identities $<\lll x,y>=<x,\bar{\lll}y>$, $<a\lll,y>=<a,y\bar{\lll}>$:
\begin{alignat*}{5}
&c'=(4|a|^2+|b|^2)c-4<a,c>a-2(a\bar{d})b+(a\bar{b})d\\ 
&d'=(4|b|^2+|a|^2)d-4<b,d>b-2(b\bar{c})a+(b\bar{a})c.
\end{alignat*}
Since 
\begin{alignat*}{5}
4<b,d>a-2(a\bar{d})b+(a\bar{b})d=&
2a(\bar{b}d)+2a(\bar{d}b)-2(a\bar{d})b+(a\bar{b})d\\
=&
2(a\bar{b})d-2[a,\bar{b},d]-2[a,\bar{d},b]+(a\bar{b})d=3(a\bar{b})d 
\end{alignat*}
and $4<a,c>b-2(b\bar{c})a+(b\bar{a})c=3(b\bar{a})c,$
we in fact have
\begin{alignat}{5}
\label{c'} c'=&(4|a|^2+|b|^2)c+3(a\bar{b})d -4g\big((a,b),(c,d)\big)a\\
\label{d'} d'=&(|a|^2+4|b|^2)d+3(b\bar{a})c-4g\big((a,b),(c,d)\big)b.
\end{alignat}
The inner product terms at the end of (\ref{c'}) and (\ref{d'}) can be ignored whenever we consider the restriction of $\jac$ to $(a,b)^\perp$.

We now describe the eigenspaces of $\jac$ corresponding to the eigenvalues $4\eps_{(a,b)}$ and $\eps_{(a,b)}$ where $\eps_{(a,b)}:=g\big((a,b),(a,b)\big)$ . In what follows we assume $|a|^2\ne 0$; we may do so because 
$|a|^2+|b|^2\ne 0$ and $(u,v)\mapsto (v,u)$ induces an isometry on  $U_1/_{\sim}$.
The solutions $(c,d)\in(a,b)^\perp$ of the linear system of equations  
$$(a\bar{b})d=|b|^2c,\ \ (b\bar{a})c=|a|^2d$$
yield eigenvectors corresponding to $4\eps_{(a,b)}$. Under our assumptions we get $$d=\frac{1}{|a|^2}(b\bar{a})c.$$ Since
$(a\bar{b})d=\frac{1}{|a|^2}(a\bar{b})(b\bar{a})c=|b|^2c$ and
\begin{alignat*}{5}
<a,c>+<b,d> 
=&<a,c>+\frac{1}{|a|^2}<b\bar{c},b\bar{a}>\\
=&<a,c>+\frac{|b|^2}{|a|^2}<a,c>=\frac{\eps_{(a,b)}}{|b|^2}<a,c>,
\end{alignat*}
we see that the $4\eps_{(a,b)}$-eigenspace is of dimension $7$ and is equal to 
$$\ker\{\jac-4\eps_{(a,b)}Id\}=\Big\{\big(c,\frac{1}{|a|^2}(b\bar{a})c\big)\Big|\ c\perp a\Big\}.$$
An analogous argument gives us
$$\ker\{\jac-\eps_{(a,b)}Id\}=\Big\{\big(-\frac{1}{|a|^2}(a\bar{b})d,d \big) \Big\}\subset (a,b)^\perp,$$
and the eigenvalue $\eps_{(a,b)}$ is of multiplicity $8$. Hence $\jac$ is diagonalizable with two non-zero eigenvalues. For the rest of the proof set $\lll=4$ and $\mu=1$.

Note that 
$$E_\lll\big((a,b)\big)=\Big\{\big(c,\frac{1}{|a|^2}(b\bar{a})c\big)\Big\}.$$
Since $\big((b\bar{a})c\big)\bar{c}=|c|^2(b\bar{a})$ and 
$|c|^2\ne 0$ for a unit vector $(c,d)\in E_\lll\big((a,b)\big)$, we have 
$$E_\lll\big((c,d)\big)=
\Big\{ \big(x,\frac{1}{|c|^2}(d\bar{c})x\big)\Big\}=
\Big\{ \big(x,\frac{1}{|a|^2}(b\bar{a})x\big)\Big\}=
E_\lll\big((a,b)\big).
$$

Finally, let $(c,d)\in \ker\{\jac-\eps_{(a,b)}Id\}$. We have $(c,d)\perp (a,b)$, $c=-\frac{1}{|a|^2}(a\bar{b})d$ and  $$|c|^2=\frac{|b|^2|d|^2}{|a|^2}.$$ 
We verify that $(a,b)\in \ker\{\J_{(c,d)}-\eps_{(c,d)}Id\}$ by a direct computation: 
\begin{alignat*}{5}
\J_{(c,d)}(a,b)=&\Big((4|c|^2+|d|^2)a+3(c\bar{d})b,\ (4|d|^2+|c|^2)b+3(d\bar{c})a\Big)\\
=&\Big((4|c|^2+|d|^2)a-3\frac{|b|^2|d|^2}{|a|^2}a,\ (4|d|^2+|c|^2)b-3|d|^2b\Big)\\
=&(|c|^2+|d|^2)(a,b)=\eps_{(c,d)}(a,b).
\end{alignat*}
\end{proof}

\begin{rem}\label{rakicsymm}
One could also verify that 
$(c,d)\in \ker\{\jac-4\eps_{(a,b)}Id\}$ implies  $(a,b)\in \ker\{\J_{(c,d)}-4\eps_{(c,d)}Id\}$; the computation is completely analogous to the one exhibited above. 
\end{rem}

We now turn to the indefinite octonionic projective plane $\oct P^{(1,1)}$. As above, we will perform all of our computations at $P_0[1,0,0]$ where the metric is given by 
$$g\big((a,b),(c,d)\big)=<a,c>-<b,d>.$$
Using the expression for  the curvature tensor at $P_0$ given at the end of section  \ref{(1,1)} we see that the Jacobi operator $\J_{(a,b)}$ (corresponding to a unit tangent vector $(a,b)$) at $P_0$ takes the form $\J_{(a,b)}(c,d)=(c',d')$ where
\begin{alignat*}{5}
c'=&(4|a|^2-|b|^2)c-3(a\bar{b})d -4g\big((a,b),(c,d)\big)a\\
d'=&(|a|^2-4|b|^2)d+3(b\bar{a})c-4g\big((a,b),(c,d)\big)b.
\end{alignat*}
\begin{thm}
$\oct P^{(1,1)}$ is a special Osserman manifold.
\end{thm}
\begin{proof}
Let $(a,b)$ be a unit tangent vector at $P_0$. An argument analogous to the one in the previous proof shows that
\begin{alignat*}{5}
\ker\{\jac-4\eps_{(a,b)}Id\}=&
\begin{cases}
\Big\{\big(c,\frac{1}{|a|^2}(b\bar{a})c\big)\Big|\ c\perp a\Big\} &\text{\ if\ } |a|\ne 0,\cr
 & \cr
\Big\{\big(\frac{1}{|b|^2}(a\bar{b})d,d\big)\Big|\ d\perp b\Big\} &\text{\ if\ } |b|\ne 0,
\end{cases}\\
\ker\{\jac-\eps_{(a,b)}Id\}=&
\begin{cases}
\Big\{\big(\frac{1}{|a|^2}(a\bar{b})d,d \big) \Big\} & \text{\ if\ } |a|\ne 0,\cr
& \cr
\Big\{\big(c,\frac{1}{|b|^2}(b\bar{a})c\big)\Big\} & \text{\ if\ } |b|\ne 0.
\end{cases}
\end{alignat*}
Distinguishing the two cases, $|a|\ne 0$ and $|b|\ne 0$, is necessary as 
$|a|^2-|b|^2\ne 0$ only guarantees $|a|\ne 0$ or $|b|\ne 0$. The rest of the proof follows along the same lines as the one above, and the details will be omitted.  \end{proof}

As mentioned earlier we complete our identification with the classical models by showing that $\oct'P^2$ is not locally isotropic. 

\begin{thm}
The split octonionic plane $\oct'P^2$ is not locally isotropic.
\end{thm}

\begin{proof}
Consider vectors $v=(i+l,0)$ and $w=(1,l)$ with base point $[1,0,0]$, where $l:=(0,1)\in \oct'$. Since $|l|^2=-1$ these vectors  are non-zero null. Suppose there exists a local isometry $\mathcal{I}$ for which 
$$\mathcal{I}[1,0,0]=[1,0,0] \text{\ \ and\ \ }
d\mathcal{I}_{[1,0,0]}v=w.$$
The vector $x:=d\mathcal{I}_{[1,0,0]}(1,0)$ must be unit spacelike. Since $\J_{(1,0)}v=4v$ we have $\J_xw=4w$, i.e. $w\in\ker\{ \J_x-4\eps_xId\}$. Remark \ref{rakicsymm} now implies $$x\in\ker\{\J_w-4\eps_wId\}=\ker\{\J_w\}.$$

Write $x=(x_1,x_2)$. It follows from $\J_wx=(3x_1-3lx_2, -3x_2+3lx_1)=0$ that 
$x_1=lx_2$. Therefore $|x_1|^2=-|x_2|^2$ and $|x|^2=|x_1|^2+|x_2|^2=0$, contradicting the fact that  $x$ is spacelike.  Therefore, there is  no local isometry $\mathcal{I}$ whose differential takes $v$ to $w$.  
\end{proof}


\end{document}